\documentclass[a4paper,oneside,reqno]{amsart}

\calclayout

\usepackage{titlesec}
\titleformat{\section}[block]
  {\normalfont\Large\bfseries}
  {\thesection}
  {1em}
  {}

\titleformat{\subsection}[block]
  {\normalfont\large\bfseries}
  {\thesubsection}
  {1em}
  {}

\usepackage{enumitem}
\usepackage{etoolbox}

\newtoggle{default}\toggletrue{default}
\newtoggle{birk}\togglefalse{birk}
\newtoggle{birk_t2}\togglefalse{birk_t2}
\newtoggle{siam}\togglefalse{siam}
\newtoggle{springer}\togglefalse{springer}

\newtoggle{endwithsymbol}\togglefalse{endwithsymbol}

\toggletrue{default}
\toggletrue{endwithsymbol}

\usepackage[utf8]{inputenc}
\usepackage[english]{babel}
\usepackage{microtype} 
\usepackage{needspace} 

\usepackage{amsmath,amssymb,amsthm}
\usepackage{thmtools} 

\usepackage{hyperref} 
\hypersetup{hidelinks}
\usepackage[nameinlink,sort]{cleveref} 
\usepackage{orcidlink}

\usepackage{mathtools} 
\usepackage{upgreek}
\usepackage{accents} 
\usepackage{stmaryrd} 

\usepackage{tikz} 
\usetikzlibrary{cd}
\usetikzlibrary{arrows.meta, positioning, calc}
\usepackage{geometry}
\usetikzlibrary{calc,patterns,shapes.geometric}
\usepackage{xcolor}

\definecolor{pointred}{RGB}{195, 35, 40}
\definecolor{softorange}{RGB}{255, 220, 175}
\definecolor{cylgray}{RGB}{180, 180, 180}
\definecolor{cylgrayback}{RGB}{140, 140, 140}
\definecolor{softgray}{RGB}{230, 230, 230}
\definecolor{conebase}{RGB}{35, 140, 70}
\definecolor{softgreen}{RGB}{210, 245, 220}
\definecolor{mainred}{RGB}{195, 35, 40}

\usepackage{dsfont} 
\usepackage{pifont} 

\usepackage{ifthen}
\setlist[itemize]{labelindent=\parindent,leftmargin=*,itemsep=3pt,topsep=5pt}
\setlist[enumerate]{label=\textup{(\roman{*})},itemsep=3pt,labelindent=0pt,leftmargin=*}

\usepackage{nicematrix}

\usepackage{graphicx} 
\usepackage{subcaption}

\usepackage{booktabs} 
\usepackage{esint} 

\usepackage{bbm}
\usepackage{float}
\usepackage{color}

\usepackage{pgfpages}
\usepackage{pgfplots}
\pgfplotsset{compat=1.13}

\newcommand{\R}{\mathbb{R}} 
\newcommand{\N}{\mathbb{N}} 
\newcommand{\C}{\mathbb{C}} 

\newcommand{\Cc}[1][\infty]{\mathrm{C}_{\mathrm{c}}\ifthenelse{\equal{#1}{}}{}{^{#1}}}
\newcommand{\Lp}[2][]{\mathrm{L}_{#2\ifthenelse{\equal{#1}{}}{}{,#1}}} 

\newcommand{\indicator}{\mathds{1}} 

\newcommand{\transposed}{^{\top}}

\DeclareMathOperator{\spt}{spt}
\DeclareMathOperator{\dom}{dom}

\DeclarePairedDelimiterX{\dset}[2]{\{}{\}}{#1\,\delimsize\vert\,\mathopen{} #2}
\DeclarePairedDelimiterX{\scprod}[2]{\langle}{\rangle}{#1,#2}

\renewcommand{\Re}{\operatorname{Re}}

\newcommand{\Htopo}{\mathrm{H}}

\theoremstyle{plain}
\newtheorem{theorem}{Theorem}[section]
\newtheorem{lemma}[theorem]{Lemma}
\newtheorem{proposition}[theorem]{Proposition}
\newtheorem{corollary}[theorem]{Corollary}
\newtheorem{hypotheses}[theorem]{Hypotheses}

\iftoggle{endwithsymbol}{%
    \declaretheorem[style=definition,sibling=theorem,qed=\ding{169}]{definition}
    \declaretheorem[style=definition,sibling=theorem,qed=\ding{169}]{example}
    \declaretheorem[style=definition,sibling=theorem,qed=\ding{169}]{problem}
    \declaretheorem[style=definition,sibling=theorem,qed=\ding{169}]{assumption}
    \declaretheorem[style=definition,numbered=no,qed=\ding{169}]{claim}

    \declaretheorem[style=remark,sibling=theorem,qed=\ding{169}]{remark}
  }%
  {%
    \theoremstyle{definition}
    \newtheorem{definition}[theorem]{Definition}
    \newtheorem{example}[theorem]{Example}

    \theoremstyle{remark}
    \newtheorem{remark}[theorem]{Remark}
  }%

\allowdisplaybreaks

\begin{document}

\title{Finite propagation speed for evolutionary equations}

\author[J.~Simo Tekeu]{Joseline Simo Tekeu}
\email{Joseline-Doriane.Simo-Tekeu@math.tu-freiberg.de}
\address{TU Bergakademie Freiberg \\
 Institute of Applied Analysis \\
Akademiestrasse 6 \\
  D-09596 Freiberg \\
  Germany}

\author[M.~Waurick]{Marcus Waurick\,\orcidlink{0000-0003-4498-3574}}
\email{marcus.waurick@math.tu-freiberg.de}

\address{TU Bergakademie Freiberg \\
  Institute of Applied Analysis \\
  Akademiestrasse 6 \\
  D-09596 Freiberg \\
  Germany}



\begin{abstract}
 We study finite speed of propagation in the framework of evolutionary equations introduced by Rainer Picard in 2009. We establish sufficient conditions on the material law, which is represented by an analytic and bounded operator-valued function, and the unbounded spatial operator that ensure finite speed of propagation of information for solutions. Assuming a linear-in-frequency coercivity bound for spatially local material laws, we prove propagation on backward light cones for weighted $L_2$ solutions without a finite initial time, thus identifying hyperbolic-type evolutionary equations. The characteristic-speed threshold depends on the coercive coefficient and the principal spatial matrices, but not on bounded memory or delay terms or potential positive fractional temporal contributions. We illustrate the result with the scalar wave equation, Maxwell's equations, a fractional model in elasticity, and dual-phase-lag thermoelasticity.    
  \end{abstract}


\ifboolexpr{togl{default} or togl{birk_t2} or togl{birk}}{\maketitle}{}%

\section{Introduction}

Finite propagation speed is a fundamental feature of hyperbolic evolution equations: a disturbance can affect another spatial point only after a positive amount of time. In the classical wave equation, Maxwell's equations and elastodynamics, this property is expressed by a domain of dependence, or equivalently by a forward or backward region of influence. It provides information beyond temporal causality. Whereas causality excludes a dependence on the future, finite propagation additionally restricts the spatial locations from which a solution can be influenced within a given time interval. For a forced equation, the appropriate statement concerns spacetime: if the forcing vanishes on a backward light cone, then so does the solution.

The distinction between finite and infinite propagation also arises in models of transport and heat conduction. Fourier's heat equation has an instantaneous spatial influence, while the telegraph-type processes studied by Goldstein~\cite{Goldstein1951} and the relaxation of Fourier's law proposed by Cattaneo~\cite{cattaneo1958forme} lead to finite-speed models. Gurtin and Pipkin~\cite{GurtinPipkin1968} developed a theory of heat conduction with memory and finite wave speeds; see also the subsequent treatment of material response functions by Finn and Wheeler~\cite{FinnWheeler1972}. Further developments include frame-indifferent formulations of Maxwell--Cattaneo heat conduction~\cite{Christov2009}, random-walk descriptions~\cite{Keller2004}, and thermoelastic models with finite wave speeds~\cite{ignaczak2010thermoelasticity}. For nonlocal models the situation is more varied: finite as well as infinite propagation can occur, depending on the constitutive and spatial operators; see, for instance, \cite{antontsev2016finite,stana2016finite,moussaoui2020speed,cheng2023dynamics,segatti2026alpha}.

Analytically, finite propagation for first-order systems is closely related to the theory of symmetric hyperbolic equations and local energy identities, going back to Friedrichs~\cite{Friedrichs1954,Friedrichs1958}. The corresponding energy or flux estimates have been adapted to systems with memory-type constitutive relations, including hyperbolic integro-differential equations and viscoelasticity~\cite{grimmer1984energy,McLaughlinYoon2017}. McIntosh and Morris~\cite{McIntoshMorris2013} obtain finite propagation for strongly continuous groups generated by suitable first-order operators. Blazek, Stolk and Symes~\cite{BlazekStolkSymes2008,BlazekStolkSymes2013} establish analogous results for first-order hyperbolic systems with bounded measurable spatial coefficients and $L_1$-convolution kernels. Finster, Murro and Schmid~\cite{FinsterMurroSchmid2025} study retarded potentials, which may also be spatially nonlocal, and formulate their finite-propagation proposition for classical solutions.

Our perspective is that of \emph{evolutionary equations}, introduced by Picard~\cite{picard2009structural} and developed systematically in \cite{PicardMcGhee2011,seifert2022evolutionary}. In exponentially weighted spacetime Hilbert spaces, the central object is the operator
\[
  \overline{\partial_{t,\nu}M(\partial_{t,\nu})+A},
\]
where the holomorphic material law $M$ describes the constitutive response, possibly with temporal nonlocality, and the operator $A$ represents the spatial differential structure and boundary conditions. This formalism treats, among others, hyperbolic, parabolic, mixed-type and integro-differential equations within a common solution theory. Its underlying research programme is to identify \emph{structural properties of $M$ and $A$ which entail properties of the corresponding solution operator}, rather than proving each property separately by constructing model-specific estimates. The basic example is Picard's coercivity criterion for well-posedness and causality. Other instances are structural criteria for maximal regularity~\cite{picard2017maximal,TrostorffWaurick2021}, exponential stability~\cite{Trostorff2015}, and continuous dependence on the coefficients, including weak-operator limits relevant to homogenisation~\cite{BuchingerSkrepekWaurick2024}. The monograph~\cite{seifert2022evolutionary} provides a broader account of these developments.

The purpose of the present paper is to add \emph{finite propagation speed} to this list of solution properties deduced from assumptions on the operator coefficients. We consider spatially local $L_\infty(\Omega)^{m\times m}$-valued material laws and skew-selfadjoint realisations of a first-order differential operator $A\subseteq \sum_{j=1}^dA_j\partial_j$ with selfadjoint matrices $A_j$. The decisive hypothesis is a quantitative, spatially pointwise coercivity bound
\begin{equation}\label{eq:intro-coercivity}
 \operatorname{Re}\bigl(zM(z)\bigr)
 \ge (\operatorname{Re}z)M_0-c_0\text{ on }L(L_2(\Omega)^m)
\end{equation}
valid for $z$ in a suitable right half-plane, where $M_0=M_0^*\in L_\infty(\Omega)^{m\times m}$ is uniformly strictly positive definite and $c_0\ge0$. Under this assumption, we prove finite propagation with any speed larger than
\[
 c_*=\operatorname*{ess\,sup}_{x\in\Omega}\sup_{|n|=1}
 \left\|M_0(x)^{-1/2}\Bigl(\sum_{j=1}^d n_jA_j\Bigr)
 M_0(x)^{-1/2}\right\|.
\]
The speed bound depends on $M_0$ and $A$, but not on the remaining temporal constitutive terms. In particular, the result covers the case $M(z)=M_0+z^{-1}M_1(z)$ with $M_1$ bounded and also material laws for which $z(M(z)-M_0)$ is unbounded, such as positive fractional time-derivative contributions. 

The proof uses Picard's well-posedness theorem together with exponential conjugation by spatial weights. A resulting off-diagonal resolvent estimate is converted, by the Fourier--Laplace transform and causality, into an exact time-delay statement for separated spatial sets. We then infer propagation on backward light cones by a countable approximation with spacetime cylinders. The conclusion applies directly to solutions in $L_{2,\nu}(\mathbb R; L_2(\Omega)^m)$ with forcing histories from the same space and does not require classical differentiability or the choice of an initial time. The result also covers holomorphic material laws with discrete delays, whose time kernels need not belong to $L_1$.

The abstract theorem is illustrated by the scalar wave equation, Maxwell's equations with variable-order fractional constitutive terms and delays, a fractional stress-compliance model in elasticity, and dual-phase-lag thermoelasticity. In the last example, we compare our characteristic-speed bound with the sufficient estimate of Mukhopadhyay, Kothari and Kumar~\cite{mukhopadhyay2011domain}. The fractional elasticity application is related to the setting of~\cite{PicardTrostorffWaurick2015fractional}. \Cref{sec:framework} recalls the evolutionary framework; \Cref{sec:finite-speed} contains the finite-propagation theorem and its proof. \Cref{sec:examples} treats the four classes of applications and compares the thermoelastic speed bounds.

\begin{section}{Framework of evolutionary equations }\label{sec:framework}
We recall the evolutionary setting and the results needed below; further details can be found in~\cite{PicardMcGhee2011,seifert2022evolutionary}. 

We begin by introducing the Hilbert space setting. 
\begin{definition}
    For $\nu\in \R$ and a Hilbert space $\Htopo$, we define the space $L_{2,\nu}(\mathbb{R}; \Htopo)$ as the space of (equivalence classes of) measurable functions $f:\R\to\Htopo$ which are square-integrable with respect to the exponentially weighted Lebesgue measure $\mathrm{e}^{-2\nu t}\mathrm{d}t,$ equipped with the inner product
\[\langle f,h \rangle_{L_{2,\nu}(\mathbb{R}; \Htopo)}\coloneqq \int_{\R}\langle f(t),h(t) \rangle_{\Htopo} \mathrm{e}^{-2\nu t}\mathrm{d}t.\]
The weight $\nu$ determines the direction of temporal causality.
Let $H_0,H_1$ be Hilbert spaces and let
$F:L_{2,\nu}(\R;H_0)\to L_{2,\nu}(\R;H_1)$.
We call $F$ \textbf{causal}, if,
for every $a\in\R$ and every $f,g\in L_{2,\nu}(\R;H_0)$ satisfying
$f=g$ almost everywhere on $(-\infty,a]$, we have
$F(f)=F(g)$ almost everywhere on $(-\infty,a]$.
\end{definition}
\begin{remark}
From the definition of $L_{2,\nu}(\R,\Htopo)$, we see 
$L_{2,0}(\R,\Htopo)= L_2(\mathbb{R}; \Htopo)$ and 
\[
\mathrm{e}^{-\nu\mathrm{m}}:L_{2,\nu}(\R;\Htopo)\to L_2(\R;\Htopo),
\quad
(\mathrm{e}^{-\nu\mathrm{m}}f)(t)=\mathrm{e}^{-\nu t}f(t),
\quad t\in\R,
\]
is unitary. Here, $\mathrm{m}$ denotes the multiplication-by-argument operator corresponding to the time parameter. 
\end{remark}

\begin{definition}
Let $\nu>0$. We define for $f \in L_{2,\nu}(\R;\Htopo)$ with $f'\in L_{2,\nu}(\R;\Htopo)$ (in the usual distributional sense), the operator $\partial_{t,\nu}f\coloneqq f'$ as an operator in $L_{2,\nu}(\R;\Htopo)$.

For $\varphi\in{C}_c^{\infty}(\R;\Htopo)$, $t\in\R$ and $\nu>0$, we recall from~\cite{picard2009structural} that 
\[
(\partial_{t,\nu}^{-1}\varphi)(t)
= \int_{-\infty}^{t} \varphi(s)\,\mathrm{d}s.\qedhere
\]
\end{definition}

\begin{remark}
As established in~\cite[Corollary~3.2.6]{seifert2022evolutionary}, the operator $\partial_{t,\nu}$ is normal and satisfies
\[
 \partial_{t,\nu}^*=-\partial_{t,\nu}+2\nu. 
\]
Furthermore, $\operatorname{Re}\partial_{t,\nu}=\nu$, and hence $\|\partial_{t,\nu}^{-1}\|\leq\frac{1}{\nu}$ for $\nu>0.$
\end{remark}

\begin{definition}[{{\cite[Section~5.3]{seifert2022evolutionary}}}]
 A mapping $M :{\dom (M)}\subseteq \mathbb{C}
\to \operatorname{L(\Htopo)}$ is called a \textbf{material law} if 
\begin{enumerate}
    \item[(a)] ${\dom }(M)$ is open and $M$ is holomorphic, and
    \item[(b)] there exists some $\nu\in\R$ such that $\mathbb{C}_{\Re>\nu}\subseteq{\dom (M)}$ and
    \[\|M\|_{\infty,\mathbb{C}_{\Re>\nu}}\coloneqq \operatorname{\sup_{z \in \mathbb{C}_{\operatorname{Re} > \nu}}}\|M(z)\|<\infty.\]
\end{enumerate}
The set
\[s_b(M)\coloneqq \inf \left\{ \nu\in \R; \ \C_{\Re>\nu}\subseteq\dom (M)\ \text{and}\ \|M\|_{\infty,\mathbb{C}_{\Re>\nu}}<\infty \right\}\]
is called the \textbf{abscissa of boundedness} of $M.$
\end{definition}

Next, as defined in~\cite[Proposition~5.3.2]{seifert2022evolutionary}, for a material law $M$ and $\nu>s_b(M)$, the \textbf{material law operator} is given by
\[
M(\partial_{t,\nu})
\coloneqq
\mathcal{L}_\nu^*
M\bigl(\mathrm{i}\mathrm{m}+\nu\bigr)
\mathcal{L}_\nu,
\] where
\[
\mathcal{L}_\nu
\coloneqq
\mathcal{F}\mathrm{e}^{-\nu\mathrm{m}}
:
L_{2,\nu}(\R;\Htopo)
\to
L_2(\R;\Htopo)
\]
and $\mathcal{F}:L_2(\R;\Htopo)\to L_2(\R;\Htopo)$ is the Fourier transformation given as the unitary extension of
\[
{C}_c^{\infty}(\R;\Htopo)
\ni
\varphi
\mapsto
\left(
s
\mapsto
\frac{1}{\sqrt{2\pi}}
\int_{\R}
\mathrm{e}^{-\mathrm{i}s t}\varphi(t)\,dt
\right).
\]
Many initial-boundary value problems can be written as evolutionary equations of the form
\[(\overline{\partial_{t,\nu}M(\partial_{t,\nu})+A})u=f,\]
where $A$ is a skew-selfadjoint operator in $\Htopo$, canonically lifted to $L_{2,\nu}(\mathbb{R};\Htopo)$.

We are now in a position to state the fundamental well-posedness theorem for evolutionary equations.
\begin{theorem}[{{\cite[Theorem~6.2.1]{seifert2022evolutionary}}}]
 \label{thm:wellposedness}
 Let $\nu_0 \in \mathbb{R}$ and $\Htopo$ be a Hilbert space. Let $M :{\dom (M)}\subseteq \mathbb{C}
\to {L(\Htopo)}$ be a material law with $s_b(M)<\nu_0$ and let $A:{\dom (A)}\subseteq  \Htopo
\to \Htopo$ be skew-selfadjoint. Assume that
\[\operatorname{Re}\langle \phi,z M(z)\phi\rangle_{\Htopo }  \ge c\left\lVert \phi\right \rVert^2_{\Htopo } \quad \left(\phi \in \Htopo, z\in \mathbb{C}_{\operatorname{Re}\ge \nu_0}\right) \]
for some $c>0.$ Then for all $\nu \ge \nu_0$ the operator $\partial_{t,\nu}M(\partial_{t,\nu})+A$ is closable and 
\[ S_{\nu} \coloneqq (\overline{\partial_{t,\nu}M(\partial_{t,\nu})+A})^{-1}\in L(L_{2,\nu}(\mathbb{R}; \Htopo)).\] Furthermore, $S_{\nu}$ is causal and satisfies $\lVert S_{\nu} \rVert_{L(L_{2,\nu})} \le \frac{1}{c},$ and for all $f\in{\dom }(\partial_{t,\nu})$ we have 
\[S_{\nu}f \in{\dom }(\partial_{t,\nu})\cap{\dom (A)}.\]
 Moreover, for $\eta, \nu\ge \nu_0$ and $f\in L_{2,\nu}(\mathbb{R}; \Htopo)\cap L_{2,\eta}(\mathbb{R}; \Htopo)$ we have that $S_{\nu}f=S_{\eta}f.$
\end{theorem}

We recall the following result on material-law operators from
\cite[Proposition~8.1.4 and Lemma~5.3.6]{seifert2022evolutionary}.

\begin{proposition}\label{prop:STW-causality}
Let $M:\dom(M)\subseteq\mathbb C\to L(H)$ be a material law.
Then for $\nu>s_b(M)$ we have
$M(\partial_{t,\nu})\in L(L_{2,\nu}(\mathbb R;H))$ and
$M(\partial_{t,\nu})$ is causal and autonomous; i.e., commutes with time-translation.
\end{proposition}

\begin{definition}\label{def:evolutionary-pair}
We call $(M,A)$ an \textbf{evolutionary pair} (in a Hilbert space $\Htopo$) if $M$ is a material law with values in $L(\Htopo)$, $A\colon \dom(A)\subseteq \Htopo \to \Htopo$ skew-selfadjoint, and there are $\nu_0>s_b(M)$ and $k>0$ such that
\[
  \operatorname{Re}\langle zM(z)v,v\rangle_{\Htopo}\ge k\|v\|_{\Htopo}^2
  \qquad(\operatorname{Re}z\ge\nu_0,\ v\in\Htopo).\qedhere
\]
\end{definition}

\begin{remark}\label{rem:resd} Let $(M,A)$ be an evolutionary pair. Then, by \cite[Theorem 6.2.1 and Remark 6.2.4]{seifert2022evolutionary}, for all $f\in L_{2,\nu}(\R;\Htopo)$, we have
\[u_{\varepsilon}\coloneqq \left(1+\varepsilon\partial_{t,\nu}\right)^{-1}u\in\dom(\partial_{t,\nu})\cap\dom(A)\text{ and }u_\varepsilon =\big(\overline{\partial_{t,\nu} M(\partial_{t,\nu})+A}\big)^{-1}  \left(1+\varepsilon\partial_{t,\nu}\right)^{-1}f,\]
where $u\coloneqq \big(\overline{\partial_{t,\nu} M(\partial_{t,\nu})+A}\big)^{-1} f$.
\end{remark}

\end{section}

\begin{section}{Finite propagation speed}\label{sec:finite-speed} 

In this section, we provide sufficient conditions to obtain finite speed of propagation of information for evolutionary equations. Throughout this section, let $\Omega\subseteq \R^d$ be open and define $\Htopo\colon L_2(\Omega)^m$ for some $d,m\in \N_{>0}$.
Before proceeding, we present a few preliminary definitions.

\begin{definition}\label{def:finite-speed}
For $x\in\Omega$ and $r\ge0$, let
\[
 B[x,r]\coloneqq\{x'\in\Omega:\ \|x-x'\|\le r\}
\]
be the closed ball in $\Omega$ with centre $x$ and radius $r$.
For $c>0$, $t_0\in\R$ and $x_0\in\Omega$, the
\textbf{backward cone} is defined by
\[
 K_c^-(t_0,x_0)\coloneqq
 \bigl\{(t,x)\in\R\times\Omega:\ t\le t_0,\ \|x-x_0\|\le c(t_0-t)\bigr\}.
\]
An evolutionary pair $(M,A)$ in $\Htopo$ has
\textbf{finite propagation speed bounded by $c_*$} (\textbf{ finite speed}, for short) if, for every $c>c_*$,
there is $\nu_0>s_b(M)$ such that, for all $\nu\ge\nu_0$,
$f\in L_{2,\nu}(\R;\Htopo)$ and $(t_0,x_0)\in\R\times\Omega$
\begin{equation}\label{eq:cone-propagation}
  \indicator_{K_c^-(t_0,x_0)}f=0
  \quad\Longrightarrow\quad
  \indicator_{K_c^-(t_0,x_0)} S_\nu f=0,
\end{equation}
where $S_\nu=\bigl(\overline{\partial_{t,\nu}M(\partial_{t,\nu})+A}\bigr)^{-1}$.
\end{definition}

\begin{remark}
    Finite speed is stronger than causality, which only yields the implication
\begin{equation}\label{eq:causality}
\indicator_{(-\infty,t_0]\times \Omega}f=0 \implies \indicator_{(-\infty,t_0]\times \Omega}\overline{(\partial_{t,\nu}M(\partial_{t,\nu})+A})^{-1}f=0.\end{equation}
The geometric distinction is illustrated in \Cref{fig:causality}.
\end{remark}
In the following, we identify $L_{\infty}$ with the corresponding multiplication operators on $L_2.$
\begin{hypotheses}\label{hyp:main}
Recall $\Htopo=L_2(\Omega)^m$ with $\Omega\subseteq\R^d$ is open. We assume:
\begin{enumerate}
\item[(i)] $A\colon\dom(A)\subseteq\Htopo\to\Htopo$ is skew-selfadjoint, and for $j\in \{1,\ldots, d\}$ there is
$A_j=A_j^*\in\C^{m\times m}$ such that
\[
 Av=\sum_{j=1}^d A_j\partial_jv
 \qquad(v\in\dom(A)).
\]
We further assume that $A$ is \textbf{\textup{local}}, meaning that $\dom(A)$ is
invariant under multiplication by scalar functions in $W^{1,\infty}(\Omega)$, i.e., for $\phi\in W^{1,\infty}(\Omega)$ and $v\in\dom(A)$,
\begin{equation}\tag{L}\label{eq:local-domain}
 \phi v\in\dom(A),\qquad
 A(\phi v)=\phi Av+\sum_{j=1}^d A_j(\partial_j\phi)v.
\end{equation}
\item[(ii)] $M_0=M_0^*\in L_\infty(\Omega)^{m\times m}$ satisfies
$M_0\ge k$ for some $k>0$.
\item[(iii)] For some $\nu_b>0$,
$M\colon\C_{\operatorname{Re}>\nu_b}\to L_\infty(\Omega)^{m\times m}$
is bounded and holomorphic. For some $c_0\ge0$ and every $z$ with $\operatorname{Re}z>\nu_b$,
\begin{equation}\label{eq:coercivity}
 \operatorname{Re}\bigl(zM(z)\bigr)
 \ge (\operatorname{Re}z)M_0-c_0\text{ on }L(L_2(\Omega)^m),
\end{equation}
where the inequality is in the sense of positive definiteness and $\operatorname{Re}B=(B+B^*)/2$.
\end{enumerate}
\end{hypotheses}
\begin{remark}
The locality assumption \eqref{eq:local-domain} is satisfied by the usual
local Dirichlet and Neumann realisations of the first-order systems under
consideration.
\end{remark}

The following is our main propagation result.
\begin{theorem}\label{thm:main}
Assume \Cref{hyp:main} and define
\begin{equation}\label{eq:characteristic-speed}
 c_*\coloneqq\operatorname*{ess\,sup}_{x\in\Omega}\ \sup_{|n|=1}
 \left\|M_0(x)^{-1/2}\left(\sum_{j=1}^d n_jA_j\right)
 M_0(x)^{-1/2}\right\|.
\end{equation}
Then $(M,A)$ is an evolutionary pair with finite propagation speed bounded by $c_*$.
\end{theorem}

\begin{figure}[htbp]
\centering
\begin{tikzpicture}[scale=1.2, >=stealth]

    \def\tInf{-3.2}
    \def\tA{2.0}          
    \def\x0{1.2}          
    \def\ry{2.4}          
    \def\rx{0.35}         

    \coordinate (Vertex) at (\tA, \x0);
    \coordinate (BaseCenter) at (\tInf, \x0);
    \coordinate (TopLeft) at (\tInf, {\x0 + \ry});
    \coordinate (BottomLeft) at (\tInf, {\x0 - \ry});

    \fill[blue!8] (\tInf, -3.0) 
        to[out=5,in=175] (\tA, -2.8) 
        -- (\tA, 4.2) 
        to[out=175,in=5] (\tInf, 4.5) 
       --cycle ;
    \draw[blue!50!black, dashed] (\tInf, -3.0) -- (\tInf, 4.5);
 
    \shade[top color=softorange!90, bottom color=conebase!70, middle color=softorange!35]
        (TopLeft) -- (Vertex) -- (BottomLeft) -- cycle;

    \begin{scope}
        \clip (TopLeft) -- (Vertex) -- (BottomLeft) -- cycle;
        \draw[step=0.35cm, conebase!35, thin] (-4, -3) grid (4, 6);
        \foreach \angle in {-15, -10, -5, 0, 5, 10, 15} {
            \draw[conebase!40, thin] (Vertex) -- ++(180+\angle:6.0);
        }
    \end{scope}
    \fill[softorange!80, opacity=0.85] (BaseCenter) ellipse ({\rx} and {\ry});
    \draw[conebase!90!black, dashed] (BaseCenter) ellipse ({\rx} and {\ry});
    \draw[conebase!50, thick] ($(BaseCenter)!0.5!(Vertex)$) ellipse ({\rx*0.5} and {\ry*0.5});

    \draw[very thick, color=black!85] (TopLeft) -- (Vertex);
    \draw[very thick, color=black!85] (BottomLeft) -- (Vertex);
    \fill[gray!30, opacity=0.75] (-2.4, 3.6) to[out=15, in=165] (-0.6, 3.8) 
                                   to[out=-45, in=90] (-0.2, 3.1)
                                   to[out=210, in=20] (-1.4, 2.9)
                                   to[out=150, in=280] (-2.4, 3.6) -- cycle;
    \draw[thick, blue!50!black] (-2.4, 3.6) to[out=15, in=165] (-0.6, 3.8) 
                                   to[out=-45, in=90] (-0.2, 3.1)
                                   to[out=210, in=20] (-1.4, 2.9)
                                   to[out=150, in=280] (-2.4, 3.6);
    \node[text=blue!60!black, font=\bfseries\large] at (-1.3, 3.35) {$f\neq0$};
    \fill[gray!30, opacity=0.75] (0.2, -1.2) to[out=20, in=170] (1.3, -1.0) 
                                 to[out=340, in=120] (1.8, -1.5)
                                 to[out=240, in=20] (1.0, -2.4)
                                 to[out=200, in=330] (0.2, -1.2) -- cycle;
    \draw[thick, blue!50!black] (0.2, -1.2) to[out=20, in=170] (1.3, -1.0) 
                                  to[out=340, in=120] (1.8, -1.5)
                                  to[out=240, in=20] (1.0, -2.4)
                                  to[out=200, in=330] (0.2, -1.2);
    \node[text=blue!60!black, font=\bfseries\large] at (1.0, -1.65) {$f\neq0$};
    \node[fill=white,fill opacity=0.85,text opacity=1,
          rounded corners=2pt,inner sep=3pt,align=center,
          font=\small] at (-0.95,0.75) {$f=0$\\$u=S_\nu f=0$};
    \draw[->, thick] (-4.2, 0) -- (3.8, 0) node[right, font=\small] {$t\in\mathbb{R}$};
    \draw[->, thick] (0, -3.5) -- (0, 5.0) node[above, font=\small] {$x\in\mathbb{R}^d$};
    \draw[dashed, color=pointred, semithick] (\tA, -3.2) -- (\tA, 4.6);
    \draw[dashed, color=pointred, semithick] (0, \x0) -- (\tA, \x0);
    \node[left, font=\small, color=pointred] at (0, \x0) {$x_0$};
    \node[below right, font=\small, color=pointred] at (\tA, 0) {$t_0$};

    \node[text=blue!70!black, font=\bfseries\small] at (1., 3.2) {$(-\infty,t_0]\times\Omega$};
   
    \node[text=pointred, font=\bfseries\small] (labelbase) at (-1.1, -1.3) {
        $B[x_0, c(t_0 - t)]$
    };
    \draw[->, thick, pointred] (labelbase.west) to[out=140, in=270] (\tInf, \x0);
    \filldraw[color=pointred] (Vertex) circle (2.0pt);
    \node[above right, font=\small] at (Vertex) {$(t_0, x_0)$};
   
\end{tikzpicture}
\caption{The yellow backward light cone
 $K_c^-(t_0,x_0)$ has opening speed $c>c_*$. For finite speed, forcing may be nonzero
 outside the cone, including at earlier times, but
 $\indicator_{K_c^-}f=0$ implies $\indicator_{K_c^-}S_\nu f=0$.
 For causality, vanishing of $f$ on the entire light-blue past slab
 implies vanishing of $u$ there; the cone conclusion needs only
 vanishing of $f$ in the yellow/green cone.}
\label{fig:causality}
 \end{figure}
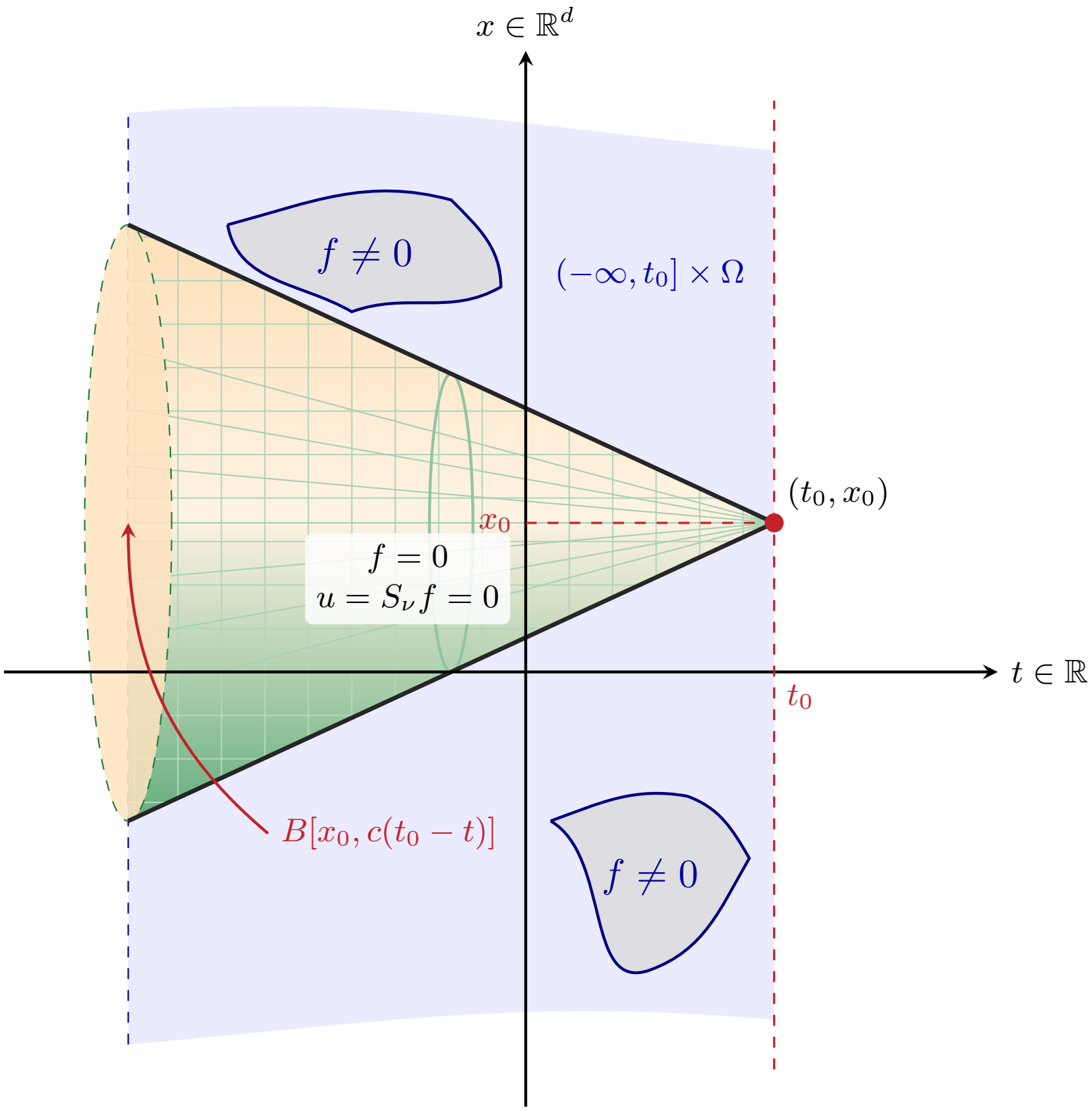

\begin{remark}\label{rem:exact-speed}
Neither the remaining part of the material law nor the constant $c_0$ in
\eqref{eq:coercivity} enters the characteristic speed $c_*$. They may
influence the admissible exponential weights, but not the propagation
speed. In particular, the finite-propagation bound is independent of
bounded memory, delay terms, and the nonnegative fractional time derivatives
covered by \Cref{hyp:main}; see the examples in \Cref{sec:examples}.
\end{remark}

The bounded-memory setting of the applications is recovered without any
positivity requirement on $M_1$.
\begin{corollary}\label{cor:bounded-memory}
Assume \Cref{hyp:main}(i)--(ii), and let
$M_1\colon\C_{\operatorname{Re}>\nu_0}\to
L_\infty(\Omega)^{m\times m}$ a material law. Define
$M(z)=M_0+z^{-1}M_1(z)$ for $\operatorname{Re}z>\max\{0,s_b(M_1)\}$.
Then $(M, A)$ has finite propagation speed bounded by $c_*$ given by \eqref{eq:characteristic-speed}.
\end{corollary}
\begin{proof}
For almost every $x\in\Omega$ and every $z$ in the given right half-plane,
\[
 \operatorname{Re}(zM(z,x))
 = (\operatorname{Re}z)M_0(x)+\operatorname{Re}M_1(z,x)
 \ge (\operatorname{Re}z)M_0(x)-\|M_1\|_\infty.
\]
On each smaller positive right half-plane $M$ is bounded and
holomorphic. Applying \Cref{thm:main} with $c_0=\|M_1\|_\infty$
proves the assertion.
\end{proof}

The distinction between bounded fractional memory and positive unbounded
fractional time derivatives is illustrated next.
\begin{remark}
For $0<\alpha<1$ and $\operatorname{Re}z>\nu_0>0$, both
\[
 M_1(z)=z^{-\alpha},\qquad
 M_1(z)=\frac{z^\alpha}{1+z^\alpha}
       =1-\frac{1}{1+z^\alpha}
\]
are bounded holomorphic material laws. In contrast,
$z^\alpha$ is unbounded on every right half-plane and thus not
covered by the corollary. However, if
$B\in L_\infty(\Omega)^{m\times m}$ is nonnegative, then
\[
 M(z,x)=M_0(x)+z^{\alpha-1}B(x)
\]
is a material law on positive right half-planes and satisfies
\eqref{eq:coercivity} with $c_0=0$, as
$\operatorname{Re}(z^\alpha)\geq (\Re z)^\alpha$; see \cite{PicardTrostorffWaurick2015fractional} or \cite[Proposition 7.2.1]{seifert2022evolutionary}.
The equation therefore admits $B\partial_{t,\nu}^\alpha$ with the
same speed bound $c_*$.
\end{remark}

We first collect the cone-locality facts used in the proof of \Cref{thm:main}.
\begin{lemma}\label{lem:invariance}
 Let $\left(\Omega,\Sigma,\mu\right)$ be a measure space, $B\in\Sigma,$ $M\in L(L_2(\mu)).$ 
 Define the closed subspace
 \[L_{2,B}(\mu)\coloneqq\left\{f\in L_{2}(\mu): \indicator_{B}f=0\right\}.\] Then,
 \[\indicator_BM=\indicator_BM\indicator_B \Longleftrightarrow M \ \text{leaves}\ L_{2,B}(\mu)\ \text{invariant}. \]
\end{lemma}
\begin{proof}
We compute
\begin{align*}
 \indicator_BM=\indicator_BM\indicator_B
 &\ \Longleftrightarrow\ \indicator_BM(1-\indicator_B)=0\\
 &\ \Longleftrightarrow\ M(\ker\indicator_B)\subseteq\ker\indicator_B\\
 &\ \Longleftrightarrow\ ML_{2,B}(\mu)\subseteq L_{2,B}(\mu).\qedhere
\end{align*}
\end{proof}
The preceding lemma will be used in the following statement.
\begin{lemma}\label{lem:union}
Let $\left(\Omega,\Sigma,\mu\right)$ be a measure space, $B,B', E\in\Sigma$, $(B_j)_{j\in \N}$ in $\Sigma$, $F\coloneqq \bigcup_{j\in \N} B_j$, and $M\in L(L_2(\mu)).$
\begin{enumerate}
\item[(a)]If $\indicator_B  M =\indicator_B  M \indicator_B  $   and $\indicator_{B'}  M =\indicator_{B'}  M \indicator_{B'},$ then
 \[\indicator_{B\cup B'}M=\indicator_{B\cup B'}M\indicator_{B\cup B'}.\]
 \item[(b)] If, for all $j\in \N$, $\indicator_{B_j}  M =\indicator_{B_j}  M \indicator_{B_j}$, then   $\indicator_{F}  M =\indicator_{F}  M \indicator_{F}$.
 \item[(c)] If, for all $j\in \N$, $\indicator_{E}  M \indicator_{B_j}=0$, then $\indicator_{E}  M \indicator_{F}=0$; if, for all $j\in \N$, $\indicator_{B_j}  M \indicator_{E}=0$, then $\indicator_{F}  M \indicator_{E}=0$.
\end{enumerate}
\end{lemma}
\begin{proof}(a) Since $\indicator_{B\cup B'}f=0$ if and only if $\indicator_Bf=0$ and $\indicator_{B'}f=0$, we have $L_{2,B\cup B'}(\mu)=L_{2,B}(\mu)\cap L_{2,B'}(\mu)$. Hence, by \Cref{lem:invariance} applied to $B$ and to $B'$,
\[
 \begin{aligned}
 M[L_{2,B\cup B'}(\mu)]
 &\subseteq M[L_{2,B}(\mu)]\cap M[L_{2,B'}(\mu)]\\
 &\subseteq L_{2,B}(\mu)\cap L_{2,B'}(\mu)
 =L_{2,B\cup B'}(\mu).
 \end{aligned}
\]
Another application of \Cref{lem:invariance}, now to $B\cup B'$, yields the first assertion.

For (b) and (c), note that for $F_N=\bigcup_{j=1}^N B_j$, $N\in \N$, we have $\indicator_{F_N}\to\indicator_{F}$ strongly by dominated convergence.

Hence, (b) follows from (a) by induction to imply $\indicator_{F_N}  M =\indicator_{F_N}  M \indicator_{F_N}$ for all $N\in \N$ and $N\to\infty$ computing the strong operator limits.

For (c), if $\indicator_E M\indicator_{B_j}=0$, we likewise obtain $\indicator_E M\indicator_{F_N}=0$ for all $N\in \N$; taking $N\to\infty$ 
proves the first assertion. For the second one, note that
$\indicator_{B_j}(M\indicator_Ef)=0$ for every $j$ means that
$M\indicator_Ef$ vanishes almost everywhere on each $B_j$, hence on
$F=\bigcup_jB_j$.
\end{proof}
The following geometric observation will be used to exhaust both the
interior and the exterior of a backward cone by countably many cylinders.
\begin{lemma}\label{lem:two-sided-cone-cylinders}
Let $t_0\in \R$, $x_0\in \Omega$, with $c>0$. Then
\[
 \mathcal{L}^{d+1}\!\left(K_c^-(t_0,x_0)\mathbin{\triangle}
 \bigcup_{\substack{s\in\mathbb Q\\s<t_0}}
 ((-\infty,s]\times B[x_0,c(t_0-s)])\right)=0.
\]
For any $s<t_0$ one also has
\[
 ((-\infty,s]\times\Omega)\setminus K_c^-(t_0,x_0)
 =\bigcup_{\substack{b\in\mathbb Q\\b<s}}
 \left((b,s]\times\{y\in\Omega:\|y-x_0\|>c(t_0-b)\}\right).
\]
\end{lemma}
\begin{proof} We abbreviate $K=K_c^-(t_0,x_0)$.
Each cylinder in the first union lies in $K$. Conversely, if
$t<t_0$ and $\|x-x_0\|<c(t_0-t)$, choose rational $s$ with
$t\le s<t_0$ and $\|x-x_0\|<c(t_0-s)$. This proves the first
assertion since $\partial K$ has Lebesgue measure zero.
For the second, if $b<t\le s$ and $\|y-x_0\|>c(t_0-b)$, then
$\|y-x_0\|>c(t_0-t)$. Conversely, if $t\le s$ and
$\|y-x_0\|>c(t_0-t)$, choose rational $b<t$ sufficiently close to
$t$ with $\|y-x_0\|>c(t_0-b)$.
\end{proof}
The first statement of the latter lemma is illustrated in \Cref{fig:cylinders}.
\begin{figure}[H]
\centering
\begin{tikzpicture}[
    scale=1.1, 
    >=stealth
]

    \draw[->, thick] (-6.5, 0) -- (2.8, 0) node[right] {$t$};
    \draw[->, thick] (-5.2, -0.5) -- (-5.2, 3.6) node[above left] {$x$};

    \coordinate (Vertex) at (1.5, 1.3);

    \def\tconeone{-4.5}   \def\ryone{1.6}    
    \def\tconetwo{-2.8}   \def\rytwo{1.15}   
    \def\tconethree{-1.2} \def\rythree{0.72} 
    \def\tcone4{0.2}      \def\ryfour{0.35}  

    \begin{scope}
        \shade[top color=purple!30, bottom color=purple!10, middle color=white, opacity=0.5] 
            (-6.5, {1.3 + \ryone}) -- (\tconeone, {1.3 + \ryone}) -- (\tconeone, {1.3 - \ryone}) -- (-6.5, {1.3 - \ryone}) -- cycle;

        \draw[purple!60!black, thick] (-6.5, {1.3 + \ryone}) -- (\tconeone, {1.3 + \ryone});
        \draw[purple!60!black, thick] (-6.5, {1.3 - \ryone}) -- (\tconeone, {1.3 - \ryone});
        
        \fill[purple!20, opacity=0.5] (\tconeone, 1.3) ellipse (0.25cm and \ryone cm);
        \draw[purple!70!black, thick] (\tconeone, 1.3) ellipse (0.25cm and \ryone cm);
    \end{scope}

    \begin{scope}
        \shade[top color=softgray!80, bottom color=gray!30, middle color=white, opacity=0.6] 
            (-6.5, {1.3 + \rytwo}) -- (\tconetwo, {1.3 + \rytwo}) -- (\tconetwo, {1.3 - \rytwo}) -- (-6.5, {1.3 - \rytwo}) -- cycle;

        \draw[gray!80, thick] (-6.5, {1.3 + \rytwo}) -- (\tconetwo, {1.3 + \rytwo});
        \draw[gray!80, thick] (-6.5, {1.3 - \rytwo}) -- (\tconetwo, {1.3 - \rytwo});
        
        \fill[gray!40, opacity=0.5] (\tconetwo, 1.3) ellipse (0.22cm and \rytwo cm);
        \draw[gray!70, thick] (\tconetwo, 1.3) ellipse (0.22cm and \rytwo cm);
    \end{scope}

    \begin{scope}
        \shade[top color=blue!40, bottom color=blue!10, middle color=white, opacity=0.6] 
            (-6.5, {1.3 + \rythree}) -- (\tconethree, {1.3 + \rythree}) -- (\tconethree, {1.3 - \rythree}) -- (-6.5, {1.3 - \rythree}) -- cycle;

        \draw[blue!60!black, thick] (-6.5, {1.3 + \rythree}) -- (\tconethree, {1.3 + \rythree});
        \draw[blue!60!black, thick] (-6.5, {1.3 - \rythree}) -- (\tconethree, {1.3 - \rythree});
        
        \fill[blue!20, opacity=0.5] (\tconethree, 1.3) ellipse (0.16cm and \rythree cm);
        \draw[blue!70!black, thick] (\tconethree, 1.3) ellipse (0.16cm and \rythree cm);
    \end{scope}

    \begin{scope}
        \shade[top color=orange!50, bottom color=orange!15, middle color=white, opacity=0.6] 
            (-6.5, {1.3 + \ryfour}) -- (\tcone4, {1.3 + \ryfour}) -- (\tcone4, {1.3 - \ryfour}) -- (-6.5, {1.3 - \ryfour}) -- cycle;

        \draw[orange!80!black, thick] (-6.5, {1.3 + \ryfour}) -- (\tcone4, {1.3 + \ryfour});
        \draw[orange!80!black, thick] (-6.5, {1.3 - \ryfour}) -- (\tcone4, {1.3 - \ryfour});
        
        \fill[orange!30, opacity=0.5] (\tcone4, 1.3) ellipse (0.10cm and \ryfour cm);
        \draw[orange!90!black, thick] (\tcone4, 1.3) ellipse (0.10cm and \ryfour cm);
    \end{scope}

    \shade[top color=softgreen!80, bottom color=conebase!60, middle color=softgreen!30, opacity=0.6]
        (\tconeone, {1.3 + \ryone}) -- (Vertex) -- (\tconeone, {1.3 - \ryone}) -- cycle;

    \fill[softgreen!80, opacity=0.6] (\tconeone, 1.3) ellipse (0.25cm and \ryone cm);
    \draw[conebase!80, thick] (\tconeone, 1.3) ellipse (0.25cm and \ryone cm);

    \draw[conebase, very thick, dashed] (\tconetwo, 1.3) ellipse (0.22cm and \rytwo cm);
    \draw[conebase!70, semithick, dashed] (\tconethree, 1.3) ellipse (0.16cm and \rythree cm);
    \draw[conebase!60, thin, dashed] (\tcone4, 1.3) ellipse (0.10cm and \ryfour cm);

    \draw[very thick, conebase!90!black] (\tconeone, {1.3 + \ryone}) -- (Vertex) -- (\tconeone, {1.3 - \ryone});

    \node[gray!80!black, font=\scriptsize] at (-5.9, -0.5) {$t \in (-\infty, t_0]$};

    \draw[dashed, mainred, semithick] (1.5, -0.5) node[below] {$t_0$} -- (1.5, 3.0);
    \draw[dashed, mainred, semithick] (-5.2, 1.3) node[left, black] {$x_0$} -- (1.5, 1.3);

    \filldraw[mainred] (Vertex) circle (2.5pt) node[above right=1pt, black, font=\small] {$(t_0, x_0)$};
\end{tikzpicture}
\caption{The backward cone as a union of past cylinders.}
\label{fig:cylinders}
\end{figure}
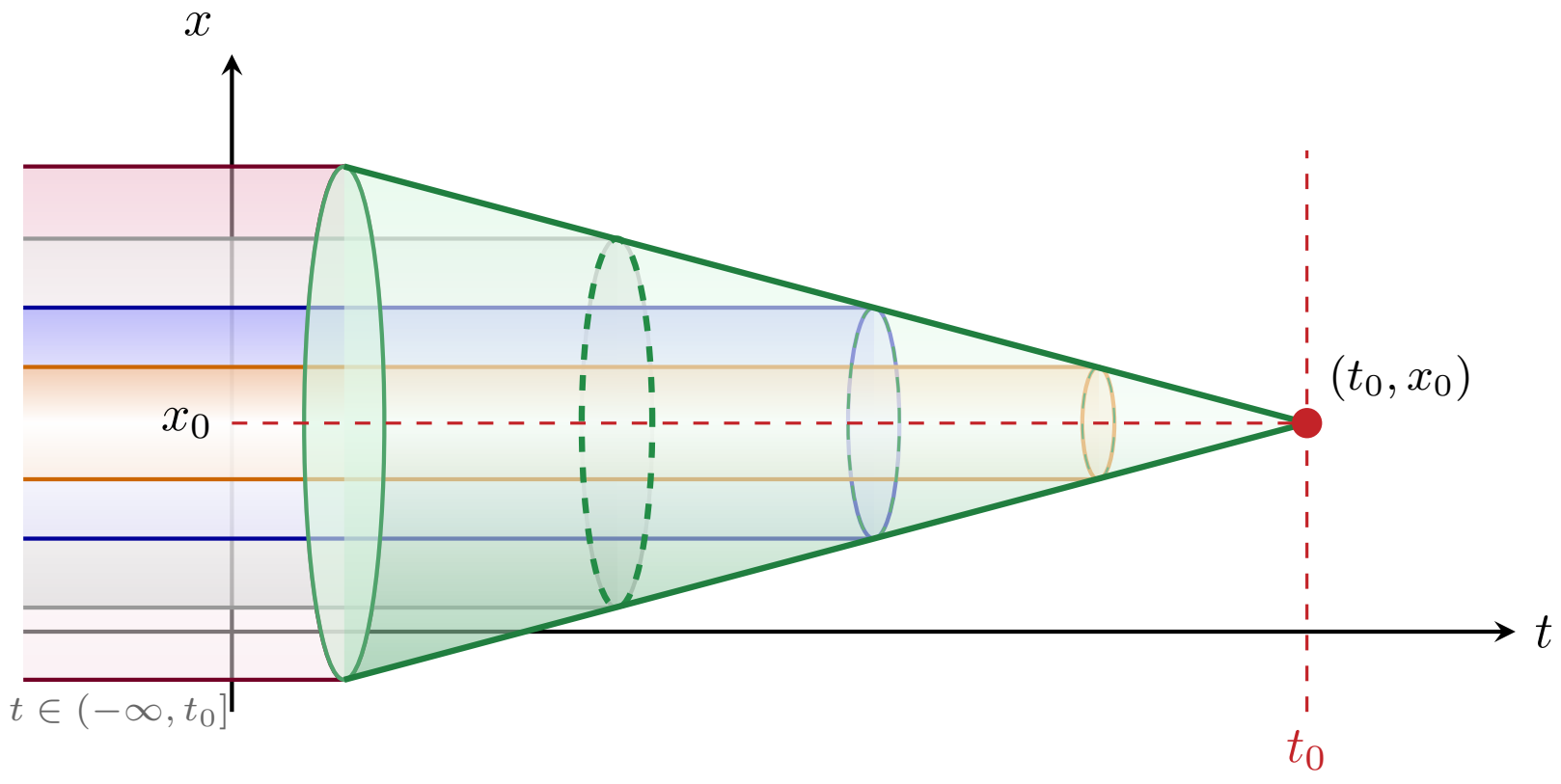
\begin{lemma}\label{lem:cone-locality}
Let $M: \mathbb{C}_{\operatorname{Re}>\nu_0}\to L_{\infty}(\Omega)^{m\times m}$
be bounded and analytic for some $\nu_0\in \mathbb{R}.$ Then for every $\nu>\nu_0$, $c>0$, $t_0\in\R$ and $x_0\in\Omega$,
\[\indicator_{K_c^-(t_0,x_0)}M(\partial_{t,\nu})=\indicator_{K_c^-(t_0,x_0)}M(\partial_{t,\nu})\indicator_{K_c^-(t_0,x_0)}.\]
\end{lemma}
\begin{proof}
Let $s\le t_0$. Since $M$ takes values in $L_\infty(\Omega)^{m\times m}$, $M(\partial_{t,\nu})$ commutes with multiplication by the scalar function $x\mapsto \indicator_{B[x_0,c(t_0-s)]}(x)$. Together with causality of $M(\partial_{t,\nu})$ (see \Cref{prop:STW-causality}), this implies
\[
\indicator_{(-\infty,s]\times B[x_0,c(t_0-s)]}M(\partial_{t,\nu})=\indicator_{(-\infty,s]\times B[x_0,c(t_0-s)]}
M(\partial_{t,\nu})
\indicator_{(-\infty,s]\times B[x_0,c(t_0-s)]}.
\]
By \Cref{lem:union}, the same identity holds for every finite union
of such past cylinders, and then for their countable union. By \Cref{lem:two-sided-cone-cylinders}, this union agrees
with $K_c^-(t_0,x_0)$ up to a null set. 
Thus
\[
\indicator_{K_c^-(t_0,x_0)}M(\partial_{t,\nu})
=\indicator_{K_c^-(t_0,x_0)}M(\partial_{t,\nu})
\indicator_{K_c^-(t_0,x_0)}.\qedhere
\]
\end{proof}

The next observation is a standard exponential-weight argument of
Combes--Thomas/Davies type (compare \cite{CombesThomas1973,CoulhonSikora2008}).
We give its short proof for reference.

\begin{lemma}\label{lem:offdiagonal}
Let $\Htopo=L_2(\Omega)^m$, $R\in L(\Htopo)$, $k>0$, $\beta>0$, and $\kappa\ge0$. Suppose that,
for every bounded real-valued $1$-Lipschitz function $\phi$ on $\Omega$
and each $0\le\eta<\beta/\kappa$ (with no upper restriction when
$\kappa=0$),
\begin{equation}\label{eq:lemma-weighted}
\|\mathrm{e}^{\eta\phi}R\mathrm{e}^{-\eta\phi}\|_{L(\Htopo)}
 \le\frac1{k(\beta-\eta\kappa)}.
\end{equation}
Then for all measurable $D,E\subseteq\Omega$ and all $c'>\kappa$,
\begin{equation}\label{eq:offdiag-abstract}
 \|\indicator_E R\indicator_D\|_{L(\Htopo)}
 \le \frac{\mathrm{e}^{-\beta\operatorname{dist}(E,D)/c'}}
 {k\beta(1-\kappa/c')}.
\end{equation}
\end{lemma}
\begin{proof}
If either set is empty, the assertion is trivial. Otherwise set $d\coloneqq \operatorname{dist}(E,D)$. 
Choose $L\ge d$ and put
$\phi_L(x)\coloneqq \min\{\operatorname{dist}(x,D),L\}$. Then $\phi_L$ is bounded and
$1$-Lipschitz, is zero on $D$, and satisfies $\phi_L\ge d$ on $E$.
Consequently,
\[
 \|\indicator_E R\indicator_D\|
 =\|\indicator_E\mathrm{e}^{-\eta\phi_L}
   (\mathrm{e}^{\eta\phi_L}R\mathrm{e}^{-\eta\phi_L})
   \mathrm{e}^{\eta\phi_L}\indicator_D\|
 \le\frac{\mathrm{e}^{-\eta d}}{k(\beta-\eta\kappa)}.
\]
Taking $\eta=\beta/c'$ gives \eqref{eq:offdiag-abstract}.
\end{proof}
We may now prove our main result.
\begin{proof}[Proof of \Cref{thm:main}]
\emph{Step 1.}
Choose an arbitrary $\nu_0$ with $\max\{\nu_b,c_0/k\}<\nu_0<\nu$. Since
\[
 \Re(zM(z))\ge (\Re z)M_0-c_0
 \ge (\Re z)k-c_0\ge k\nu_0-c_0
 \qquad (\Re z\ge\nu_0),
\]
and since $M$ is a material law and $A$ is skew-selfadjoint,
\Cref{thm:wellposedness} yields the continuous invertibility of
$\overline{\partial_{t,\nu}M(\partial_{t,\nu})+A}$,
the causality of its inverse $S_\nu$, and the estimate
$\|S_\nu\|\le(k\nu_0-c_0)^{-1}$.
By \cite[Remark~6.3.4, p.~99]{seifert2022evolutionary},
the inverse is a material-law operator. More explicitly, with
\[
  R(z)\coloneqq (zM(z)+A)^{-1},\qquad \Re z>\max\{\nu_b,c_0/k\},
\]
we have
\begin{equation}\label{eq:representation}
 S_\nu=R(\partial_{t,\nu})
 =\mathcal L_\nu^*R(\mathrm{i}\mathrm{m}+\nu)\mathcal L_\nu.
\end{equation}
The holomorphic function $R$ satisfies
\begin{equation}\label{eq:Rbound}
 \|R(z)\|\le\frac1{k\Re z-c_0}
 \qquad(\Re z>\max\{\nu_b,c_0/k\}).
\end{equation}
Indeed, for $u=R(z)f$, we have
$(k\Re z-c_0)\|u\|^2\le\Re\langle(zM(z)+A)u,u\rangle\le\|f\|\|u\|$.
\par\medskip\noindent
\emph{Step 2.}
Let $\phi:\Omega\to\R$ be bounded and $1$-Lipschitz,
and $\eta\ge0$.
Since $\mathrm{e}^{\pm\eta\phi}\in W^{1,\infty}(\Omega)$,
locality of $A$ and commutation of $zM(z)$ with spatial multiplication imply
\begin{equation}\label{eq:conjugate}
 \mathrm{e}^{\eta\phi}(zM(z)+A)\mathrm{e}^{-\eta\phi}
 =zM(z)+A-\eta B_\phi,\text{ where }
 \qquad B_\phi(x)\coloneqq \sum_{j=1}^d A_j\partial_j\phi(x),
\end{equation}
as an operator identity on $\dom(A)$.
For each $\xi\in\C^m$, we estimate
\begin{equation}\label{eq:matrix}
 |\langle B_\phi(x)\xi,\xi\rangle|
 \le c_*|\nabla\phi(x)|\langle M_0(x)\xi,\xi\rangle
 \le c_*\langle M_0(x)\xi,\xi\rangle.
\end{equation}
Indeed, for $\nabla\phi(x)\ne0$ this follows from the definition of $c_*$
with $n=\nabla\phi(x)/|\nabla\phi(x)|$; for $\nabla\phi(x)=0$
it is immediate. Thus, whenever $0\le\eta c_*<\Re z-c_0/k$,
\begin{equation}\label{eq:coercivity-conjugate}
 \begin{aligned}
 \Re\langle
 \mathrm{e}^{\eta\phi}(zM(z)+A)\mathrm{e}^{-\eta\phi}w,w\rangle
 &\ge(\Re z-\eta c_*)\langle M_0w,w\rangle-c_0\|w\|^2\\
 &\ge k(\Re z-c_0/k-\eta c_*)\|w\|^2
 \qquad(w\in\dom(A)).
 \end{aligned}
\end{equation}
Conjugation preserves invertibility, and hence
\begin{equation}\label{eq:weighted-resolvent}
 \|\mathrm{e}^{\eta\phi}R(z)\mathrm{e}^{-\eta\phi}\|
 \le\frac{1}{k(\Re z-c_0/k-\eta c_*)}.
\end{equation}
For any $c'>c_*$, Lemma~\ref{lem:offdiagonal} applies to
$R=R(z)$, $\beta=\Re z-c_0/k$, and $\kappa=c_*$, in view of
\eqref{eq:weighted-resolvent}. Thus for measurable $D,E\subseteq\Omega$,
\begin{equation}\label{eq:offdiag}
 \|\indicator_E R(z)\indicator_D\|
 \le \frac{\mathrm{e}^{-(\Re z-c_0/k)\operatorname{dist}(E,D)/c'}}
 {k(\Re z-c_0/k)(1-c_*/c')}\qquad(\Re z>\max\{\nu_b,c_0/k\}).
\end{equation}
\emph{Step 3.}
We apply Proposition~\ref{prop:STW-causality} next. For this, let $D,E\subseteq\Omega$ be measurable,
fix $c'>c_*$, and put
\[
 d\coloneqq\operatorname{dist}(E,D),\qquad
 a\coloneqq\frac{d}{c'},\qquad
 R_{E,D}(z)\coloneqq\mathrm{e}^{az}\indicator_E R(z)\indicator_D.
\]
By \eqref{eq:offdiag}, for $\Re z>\max\{\nu_b,c_0/k\}$,
\[
 \|R_{E,D}(z)\|\le
 \frac{\mathrm{e}^{c_0d/(kc')}}{k(\Re z-c_0/k)(1-c_*/c')}.
\]
For every $\rho>\max\{\nu_b,c_0/k\}$ this is uniformly bounded
on $\Re z>\rho$; hence $R_{E,D}$ is a material law with
$s_b(R_{E,D})\le\max\{\nu_b,c_0/k\}$. By
Proposition~\ref{prop:STW-causality},
$R_{E,D}(\partial_{t,\nu})$ is causal for each
$\nu>\max\{\nu_b,c_0/k\}$.
Let $\tau_a g(t)=g(t-a)$. Since $\mathrm{e}^{-a\partial_{t,\nu}}=\tau_{a}$, which is easily seen by Fourier--Laplace transformation, we have
\begin{equation}\label{eq:time-delay}
 \indicator_E S_\nu\indicator_D
 =\tau_a R_{E,D}(\partial_{t,\nu}),
 \qquad a=\frac{\operatorname{dist}(E,D)}{c'}.
\end{equation}
Therefore, for any $g$ supported in $[s,\infty)\times D$,
\begin{equation}\label{eq:delay-support}
 \indicator_E S_\nu g(t)=0
 \quad\text{for almost every }t<s+\frac{\operatorname{dist}(E,D)}{c'}.
\end{equation}
\emph{Step 4.}
Let $c_*<c'<c$.
Write $K=K_c^-(t_0,x_0)$ and let $f\in L_{2,\nu}(\R;\Htopo)$
with $\indicator_K f=0$. Put $u=S_\nu f$.

Take $s<t_0$, and define $ E_s\coloneqq B[x_0,c(t_0-s)]$. We prove $\indicator_{(-\infty,s]\times E_s }u=0$.
For each $b<s$, let
\[
 D_b\coloneqq\{y\in\Omega:|y-x_0|>c(t_0-b)\},\qquad
 Q_b\coloneqq (b,s]\times D_b.
\]
By the triangle inequality,
\begin{equation}\label{eq:directgap}
  \operatorname{dist}(E_s,D_b)\ge c(s-b).
\end{equation}
For every $g$ with $\spt g \subseteq Q_b\subseteq [b,\infty)\times D_b$ we may apply \eqref{eq:delay-support}, with $E=E_s$ and $D=D_b$,  to show that $\indicator_{E_s}S_\nu g(t)=0$ for almost all
\[
  t<b+\frac{\operatorname{dist}(E_s,D_b)}{c'}.
\]
Since \[  b+\frac{\operatorname{dist}(E_s,D_b)}{c'}
  \ge b+\frac{c(s-b)}{c'}>s,\] we obtain, as $g$ was chosen arbitrary satisfying the support constraint, 
\begin{equation}\label{eq:direct-cyl-vanishing}
  \indicator_{(-\infty,s]\times E_s }S_\nu\indicator_{Q_b}=0
  \qquad (b\in\mathbb Q,\ b<s).
\end{equation}
Recall from \Cref{lem:two-sided-cone-cylinders} that
\begin{equation}\label{eq:direct-source-cover}
  ((-\infty,s]\times\Omega)\setminus K
  =\bigcup_{\substack{b\in\mathbb Q\\b<s}} Q_b.
\end{equation}
Hence, \Cref{lem:union},
\eqref{eq:direct-cyl-vanishing}, and \eqref{eq:direct-source-cover} imply
\[
\indicator_{(-\infty,s]\times E_s }S_\nu
\indicator_{(( -\infty,s]\times\Omega)\setminus K}=0.
\]
Since $S_\nu$ is causal,
$\indicator_{(-\infty,s]\times E_s }S_\nu\indicator_{(s,\infty)\times\Omega}=0$. As $f$ vanishes on $K$, we get for $u=S_\nu f$ and $s<t_0$
\[
     \indicator_{(-\infty,s]\times E_s } u = 0,
\]
as indicated above. Rewriting this, using that $f$ is not supported on $K$ and that $f$ was arbitrary, we get
\begin{equation}\label{eq:direct-target-cylinder}
  \indicator_{(-\infty,s]\times E_s }S_\nu(1-\indicator_K)=0\quad(s<t_0).
\end{equation}
Finally, \Cref{lem:two-sided-cone-cylinders} and \Cref{lem:union} applied to \eqref{eq:direct-target-cylinder} yield
\[
  \indicator_KS_\nu(1-\indicator_K)=0,
  \text{ or, equivalently, }
  \indicator_KS_\nu=\indicator_KS_\nu\indicator_K.\qedhere
\]
\end{proof}
The following observations distinguish bounded perturbations from genuinely
unbounded temporal terms covered by the main theorem.
\begin{remark}\label{rem:unbounded-fractional}
\begin{enumerate}[label=(\alph*)]
\item For the bounded memory terms of \Cref{cor:bounded-memory}, the
memory-free inverse has norm at most $(k\nu)^{-1}$. Thus, for sufficiently large $\nu$, \((\partial_{t,\nu}M_0+A)^{-1}M_1(\partial_{t,\nu})\) has norm less than one, and the corresponding Neumann series converges. This gives an alternative, perturbative view of the bounded-memory case independently of the above proposed proof of \Cref{cor:bounded-memory} once the respective (classical) result for $(\partial_{t,\nu}M_0+A)^{-1}$ is shown.
\item For the scalar transport example on $\Omega=\R$ with $M_0=B=1$,
$A=\partial_x$ and $M(z)=1+z^{\alpha-1}$, $0<\alpha<1$,
we obtain finite propagation with every $c>1$ for
\[
 (\partial_{t,\nu}+\partial_x+\partial_{t,\nu}^{\alpha})u=f.
\]
On the other hand, the natural memory-free perturbation product
$(\partial_{t,\nu}+\partial_x)^{-1}\partial_{t,\nu}^{\alpha}$
has Fourier multiplier
\[
 \frac{(\nu+\mathrm{i}\omega)^\alpha}{\nu+\mathrm{i}(\omega+\xi)}.
\]
Taking $\xi=-\omega$ shows that this multiplier is unbounded
for every $\nu>0$. Thus finite propagation in the present generality
does not follow from the corresponding operator-norm Neumann series.\qedhere
\end{enumerate}
\end{remark}
As a special case, cone invariance yields finite propagation for classical (hyperbolic, abstract) Cauchy problems.
\begin{theorem}\label{thm:ivp-finite-speed}
Let $A$ satisfy \Cref{hyp:main}(i), and let $M_1\in L_\infty(\Omega)^{m\times m}$ and $u_0\in\dom(A)$. Let $u\in C^1([0,\infty);\Htopo)$ with $u(t)\in\dom(A)$ for all $t\ge0$ be the classical solution of
\[
 \partial_tu=-(M_1+A)u,\qquad u(0)=u_0.
\]
For $t_0>0$, $x_0\in\Omega$ and any
$c>\sup_{|n|=1}\|\sum_{j=1}^d n_jA_j\|$, if
\[\indicator_{B[x_0,ct_0]}u_0=0,\]
then
\[
 u(t,x)=0\quad\text{for almost every }(t,x)\text{ with }
 0\le t\le t_0,\ \|x-x_0\|\le c(t_0-t).
\]
\end{theorem}
\begin{proof}
Extend $u$ by zero to $t<0$ and set
\[
 w\coloneqq u-\indicator_{[0,\infty)}u_0,\qquad
 f\coloneqq-\indicator_{[0,\infty)}(M_1+A)u_0.
\]
By the evolutionary initial-value formulation
(\cite[Theorem~9.4.2]{seifert2022evolutionary}),
\[
 (\overline{\partial_{t,\nu}+M_1+A})w=f.
\]
By locality of $A$ and $M_1$,
\[
 \indicator_{B[x_0,ct_0]}u_0=0
 \quad\Longrightarrow\quad
 \indicator_{B[x_0,ct_0]}(M_1+A)u_0=0.
\]
Since
\[
 K_c^-(t_0,x_0)\cap([0,\infty)\times\Omega)
 \subseteq[0,t_0]\times B[x_0,ct_0],
\]
we obtain, by \Cref{cor:bounded-memory} (applied with $M_0=1$ and $M_1(z)=M_1$), as $\indicator_{K_c^-(t_0,x_0)}f=0$ that
\[ 
 \indicator_{K_c^-(t_0,x_0)}w=0.
\]
Thus, from
\[
 \indicator_{K_c^-(t_0,x_0)}\indicator_{[0,\infty)}u_0=0,
\]we get $\indicator_{K_c^-(t_0,x_0)}u=0$ proving the claim.
\end{proof}
\end{section}
\begin{section}{Examples}\label{sec:examples} 
In this section, we illustrate the result proved in the previous section with different classes of differential equations.

\begin{subsection}{Scalar Wave Equation}
The first result concerning finite speed deals with the scalar wave equation in an open subset $\Omega\subseteq\mathbb{R}^d$: For a given external force $f \colon \mathbb{R}\times\Omega\to\mathbb{R}$, this equation reads
\[ \partial_t^2u-\operatorname{div}\sigma=f,
\]
where $u:\mathbb R\times\Omega\to\mathbb R$ and
$\sigma:\mathbb R\times\Omega\to\mathbb R^d$ are the unknowns.
We use the divergence operator
\[
\begin{array}{rcl}
\operatorname{div}\colon \Htopo(\operatorname{div},\Omega) \subseteq L_{2}(\Omega)^{d}
&\longrightarrow&
L_{2}(\Omega)\\[2mm]
(\Phi_j)_{j\in \{1,\dots,d\}}
&\longmapsto&
\displaystyle
\sum_{j=1}^{d}\partial_j\Phi_{j},
\end{array}
\] where, as usual, \[\Htopo(\operatorname{div},\Omega)\coloneqq \left\{\Phi\in L_{2}(\Omega)^{d}: \sum_{j=1}^{d}\partial_j\Phi_{j} \in L_{2}(\Omega)\text{ in the distributional sense} \right\}.\]
Let $T\in L_\infty(\Omega)^{d\times d}$ satisfy 
\[\sigma=T\operatorname{grad}u,\qquad T(x)=T(x)^*\ge c_T \quad\text{for a.e. }x\in\Omega,\]
where $c_T>0$.
We introduce
\[
\begin{array}{rcl}
\operatorname{grad}\colon H^1(\Omega)\subseteq L_2(\Omega)
&\longrightarrow&
L_{2}(\Omega)^{d}\\[2mm]
\psi
&\longmapsto&
\displaystyle (\partial_j \psi)_{j\in \{1,\dots,d\}}
\end{array}
\]
and set $v\coloneqq\partial_tu$.
We impose homogeneous Dirichlet conditions on the velocity $v$. From \cite[Introduction]{seifert2022evolutionary}, a convenient reformulation of the wave equation in evolutionary equation form is the block operator matrix equation reading as 
\[
\partial_{t,\nu}\begin{pmatrix}
1& 0  \\
0 & T^{-1} \end{pmatrix}\begin{pmatrix}
v  \\
\sigma \end{pmatrix}
+\begin{pmatrix}
0 & -\operatorname{div} \\
-\operatorname{grad_0} & 0 
\end{pmatrix}\begin{pmatrix}
v  \\
\sigma \end{pmatrix}=\begin{pmatrix}
f  \\
0 \end{pmatrix},
\]
where $\operatorname{grad_0}=-\operatorname{div^*}$  
with
\begin{equation}\label{eq:wave-system}
M(\partial_{t,\nu})=M_0=\begin{pmatrix}
1& 0  \\
0 & T^{-1} \end{pmatrix}
\ \text{and} \
A=\begin{pmatrix}
0 & -\operatorname{div} \\
-\operatorname{grad_0} & 0 
\end{pmatrix}.
\end{equation}

\begin{theorem}
Let $\Htopo=L_2(\Omega)\times L_2(\Omega)^d$ and let $(M,A)$ be as in \eqref{eq:wave-system}. Then $(M,A)$ has finite speed; more precisely, \eqref{eq:cone-propagation} holds for every
\[
 c>\operatorname*{ess\,sup}_{x\in\Omega}\|T(x)\|^{1/2}.
\]
\end{theorem}

\begin{proof}
We begin by checking the assumptions of \Cref{cor:bounded-memory}; here $M_1=0$. Next,
\begin{enumerate}[label=(\alph*)]
\item $M_0$ is self-adjoint and positive definite.
Indeed, let $(x,y) \in \operatorname{L}_2(\Omega)\times \operatorname{L}_2(\Omega)^d.$ 
\begin{align*}
\left\langle M_0 (x,y), (x,y)\right\rangle_{\operatorname{L}_2(\Omega)\times \operatorname{L}_2(\Omega)^d}
&=\left\langle x, x\right\rangle_ {\operatorname{L}_2(\Omega)}+\left\langle T^{-1}y, y\right\rangle_{\operatorname{L}_2(\Omega)^d}\\ &\ge \min \left\{ 1, \frac{c_T}{\|T\|_\infty^2} \right \} \left\lVert (x,y) \right\rVert^2_{\operatorname{L}_2(\Omega)\times \operatorname{L}_2(\Omega)^d},
\end{align*}
where we used \cite[Proposition~6.2.3(b)]{seifert2022evolutionary} in the last step.
\item The operators $\operatorname{grad}_0$ and $\operatorname{div}$ are
closed and densely defined by
\cite[Proposition~6.1.1]{seifert2022evolutionary}. Using the adjoint
formula for off-diagonal block operator matrices in
\cite[Proposition~6.2.3(a)]{seifert2022evolutionary}, we obtain
\[
 A^*=\begin{pmatrix}
 0&-\operatorname{grad}_0^*\\
 -\operatorname{div}^*&0
 \end{pmatrix}
 =\begin{pmatrix}
 0&\operatorname{div}\\
 \operatorname{grad}_0&0
 \end{pmatrix}=-A.
\]
\item $A_j$ is self-adjoint.
Indeed,
\[A=\begin{pmatrix}
0 & -\operatorname{div} \\
-\operatorname{grad_0} & 0 
\end{pmatrix}=\begin{pmatrix}
0 &  -\sum_{j=1}^{d}\partial_j d_j  \\
 -\sum_{j=1}^{d}\partial_j g_j  & 0
\end{pmatrix}\\
=\sum_{j=1}^{d} \partial_jA_j,
\]
where 
\[
\begin{array}{rcl}
g_j:\R &\longrightarrow& \R^d\\
v&\longmapsto& ve_j
\end{array},
\
\begin{array}{rcl}
d_j:\R^d &\longrightarrow& \R\\
\sigma&\longmapsto& e_j^{\top}\sigma
\end{array},
\
 A_j=-\begin{pmatrix}  0 &  d_j  \\
 g_j  & 0 
\end{pmatrix}
\] and $(e_1,\ldots,e_d)$ the canonical basis of  $\R^d.$
\end{enumerate} 
Moreover, multiplication by any scalar $\phi\in W^{1,\infty}(\Omega)$ preserves both $H_0^1(\Omega)$ and $H(\operatorname{div},\Omega)$, with the usual product rules. Thus $A$ also satisfies the local-domain condition~\eqref{eq:local-domain}.

Finally, for $n\in\R^d$ with $|n|=1$, we have $\sum_{j=1}^dn_jA_j=-\begin{pmatrix}0&n\transposed\\ n&0\end{pmatrix}$ and $M_0(x)^{-1/2}=\begin{pmatrix}1&0\\0&T(x)^{1/2}\end{pmatrix}$, so that
\[
 \Bigl\|M_0(x)^{-1/2}\Bigl(\sum_{j=1}^dn_jA_j\Bigr)M_0(x)^{-1/2}\Bigr\|
 =\left\|\begin{pmatrix}0&(T(x)^{1/2}n)\transposed\\ T(x)^{1/2}n&0\end{pmatrix}\right\|
 =\bigl|T(x)^{1/2}n\bigr|,
\]
whose supremum over $|n|=1$ equals $\|T(x)\|^{1/2}$. Hence
\Cref{cor:bounded-memory} yields the assertion.
\end{proof}

\begin{remark}  
For $d=3,$ the corresponding matrices are given as follows:
\[A_1=-\begin{pmatrix}
0 & 1&0 &0\\
1 & 0&0&0 \\
0 & 0&0&0\\
0&0&0&0
\end{pmatrix},
\
A_2=-\begin{pmatrix}
0 & 0&1 &0\\
0 & 0&0&0 \\
1 & 0&0&0\\
0&0&0&0
\end{pmatrix},
\
A_3=-\begin{pmatrix}
0 & 0&0 &1\\
0 & 0&0&0 \\
0 & 0&0&0\\
1&0&0&0
\end{pmatrix}.
\]
 The matrices $A_j,\ j\in\{1,2,3\}$, are self-adjoint.    
\end{remark}

\end{subsection}

\begin{subsection}{Maxwell's equations}

The finite speed of propagation of solutions to Maxwell's equations has been addressed in~\cite[Example~2]{grimmer1984energy} using energy estimates. We consider Maxwell's equations 
\[ \begin{cases}
\partial_t\varepsilon E+\sigma E-\operatorname{curl}H=j,\\
\partial_t \mu H+\operatorname{curl_0}E=0.
\end{cases}\]
We employ
\[
\begin{array}{rcl}
\operatorname{curl}\colon \Htopo(\operatorname{curl},\Omega)  \subseteq L_{2}(\Omega)^{3}
&\longrightarrow&
L_{2}(\Omega)^3\\[2mm]
(\Phi)_{j\in \{1,2,3\}}
&\longmapsto&
\displaystyle
\begin{pmatrix}
    \partial_2\Phi_3-\partial_3\Phi_2\\
    \partial_3\Phi_1-\partial_1\Phi_3\\
    \partial_1\Phi_2-\partial_2\Phi_1
\end{pmatrix},
\end{array}
\] where \[\Htopo(\operatorname{curl},\Omega)\coloneqq \left\{\Phi\in L_{2}(\Omega)^{3}:  \operatorname{curl} \Phi \in L_{2}(\Omega)^{3} \text{ distributionally}\right\}.\]
Furthermore, $\operatorname{curl^*=curl_0.}$

The corresponding block operator matrix representation is
\[
\partial_{t,\nu}\begin{pmatrix}
\varepsilon& 0  \\
0 & \mu \end{pmatrix}\begin{pmatrix}
E  \\
H \end{pmatrix}
+\begin{pmatrix}
\sigma & 0\\
0 & 0 
\end{pmatrix}\begin{pmatrix}
E  \\
H \end{pmatrix}
+\begin{pmatrix}
0 & -\operatorname{curl} \\
\operatorname{curl_0} & 0 
\end{pmatrix}\begin{pmatrix}
E  \\
H \end{pmatrix}=\begin{pmatrix}
j \\
0 \end{pmatrix},
\]
with
\begin{equation}\label{eq:maxwell-system}
M(\partial_{t,\nu})=\begin{pmatrix}
\varepsilon& 0  \\
0 &  \mu \end{pmatrix}
+\partial_{t,\nu}^{-1}\begin{pmatrix}
\sigma& 0  \\
0 & 0 \end{pmatrix}=M_0+\partial_{t,\nu}^{-1}M_1, 
\
A=\begin{pmatrix}
0 & -\operatorname{curl} \\
\operatorname{curl_0} & 0 
\end{pmatrix},
\end{equation}
 $\varepsilon =\varepsilon^{*}\in  L_{\infty}(\Omega)^{3\times 3},$ $\mu=\mu^* \in  L_{\infty}(\Omega)^{3\times 3}$, $\sigma\in L_\infty(\Omega)^{3\times3}$, and $\varepsilon\ge c_1,$ $\mu\ge c_2$ for some $c_1,c_2>0$.
The matrices representing the curl operator are
\[C_1=\begin{pmatrix}
0 & 0&0 \\
0 & 0&-1 \\
0 & 1&0
\end{pmatrix},
\
C_2=\begin{pmatrix}
0 & 0&1 \\
0 & 0&0 \\
-1 & 0&0
\end{pmatrix},
\ \text{ and }
C_3=\begin{pmatrix}
0 & -1&0 \\
1 & 0&0 \\
0 & 0&0
\end{pmatrix}.
\]
\begin{theorem}
Let $\Htopo=L_2(\Omega)^3\times L_2(\Omega)^3$ and let $(M,A)$ be as in \eqref{eq:maxwell-system}. Then $(M,A)$ has finite speed; more precisely, \eqref{eq:cone-propagation} holds for every
\[
 c>\operatorname*{ess\,sup}_{x\in\Omega}\ \sup_{|n|=1}\Bigl\|\varepsilon(x)^{-1/2}\Bigl(\sum_{j=1}^3n_jC_j\Bigr)\mu(x)^{-1/2}\Bigr\|,
\]
In particular, $c>(c_1c_2)^{-1/2}$ is sufficient.
\end{theorem}
\begin{proof}
Again, we verify the assumptions of \Cref{cor:bounded-memory}. It holds that
\begin{enumerate}[label=(\alph*)]
 \item $M_0$ is self-adjoint and positive definite. Indeed, let $(\phi,\psi) \in \operatorname{L}_2(\Omega)^3\times \operatorname{L}_2(\Omega)^3.$ 
\begin{equation*}
\begin{aligned}
\left\langle M_0 (\phi,\psi), (\phi,\psi) \right\rangle_{\operatorname{L}_2(\Omega)^3\times \operatorname{L}_2(\Omega)^3}&=\left\langle \varepsilon \phi, \phi \right\rangle_ {\operatorname{L}_2(\Omega)^3}+\left\langle \mu \psi, \psi\right\rangle_{\operatorname{L}_2(\Omega)^3}\\
&\ge \min \left\{ c_1, c_2 \right \} \left\lVert (\phi,\psi) \right\rVert^2_{\operatorname{L}_2(\Omega)^3\times \operatorname{L}_2(\Omega)^3}. 
\end{aligned}
\end{equation*}
\item $A$ is skew-selfadjoint, by $\operatorname{curl}^*=\operatorname{curl}_0$ and \cite[Proposition~6.2.3(a)]{seifert2022evolutionary}.
\item $A=\sum_{j=1}^{3} \partial_j  A_j$ and $A_j$ is self-adjoint.
    
Indeed, 
\[A=\begin{pmatrix}
0 & -\operatorname{curl} \\
\operatorname{curl_0} & 0 
\end{pmatrix}=\begin{pmatrix}
0 & -\sum_{j=1}^{3}\partial_j  C_j \\
\sum_{j=1}^{3}\partial_j C_j  & 0 
\end{pmatrix}=\sum_{j=1}^{3} \partial_j  A_j 
\]
and
\[ A_j=\begin{pmatrix}  0 &  -C_j  \\
C_j  & 0 
\end{pmatrix},\ \text{ $A_j$ self-adjoint, $j\in \{1,2,3\}$.}\]
\end{enumerate}
Multiplication by any scalar $\phi\in W^{1,\infty}(\Omega)$ preserves both $H(\operatorname{curl},\Omega)$ and the homogeneous tangential-trace domain of $\operatorname{curl}_0$. The curl product rule therefore gives~\eqref{eq:local-domain}.

Finally, let $n\in\R^3$ with $|n|=1$ and set $C_n\coloneqq\sum_{j=1}^3n_jC_j$ and $B\coloneqq\varepsilon^{-1/2}C_n\mu^{-1/2}$. Since $C_n\transposed=-C_n$, we obtain
\[
 M_0^{-1/2}\Bigl(\sum_{j=1}^3n_jA_j\Bigr)M_0^{-1/2}=\begin{pmatrix}0&-B\\-B^*&0\end{pmatrix},
\]
whose norm equals $\|B\|\le(c_1c_2)^{-1/2}$, as $\|C_n\|=|n|=1$. Thus \Cref{cor:bounded-memory} yields the assertion.
\end{proof}
We consider Maxwell's equations with spatially varying fractional orders, next.
\begin{example}\label{ex:fractional-maxwell}
Consider the operator $A$ from \eqref{eq:maxwell-system} with scalar
$\varepsilon,\mu\in L_\infty(\Omega)$ satisfying
$\varepsilon,\mu\ge k>0$ for some $k>0$. Take
$\sigma,a\in L_\infty(\Omega)$ and a measurable function
$\alpha\colon\Omega\to[\alpha_-,\alpha_+]$ with
$0<\alpha_-\le\alpha_+<1$. We define
\[
 M_0(\cdot)=\begin{pmatrix}\varepsilon&0\\0&\mu\end{pmatrix} \in L(L_2(\Omega)^6)
\text{ and }
 M_1(z)(x)=\begin{pmatrix}
 \left(\sigma(x)+a(x)\dfrac{z^{\alpha(x)}}{1+z^{\alpha(x)}}\right)&0\\
 0&0
 \end{pmatrix}\quad (\operatorname{Re}z>0)
\]
using the principal branch of the fractional power. With
$w=(1+\partial_{t,\nu}^{\alpha(\cdot)})^{-1}E$, the corresponding
Maxwell's equations are
\[
\begin{cases}
 \varepsilon\partial_{t,\nu}E+\sigma E
       +a\partial_{t,\nu}^{\alpha(\cdot)}w-\operatorname{curl}H=j,\\
 \mu\partial_{t,\nu}H+\operatorname{curl}_0E=0.
\end{cases}
\]
Here fractional powers act in time with the spatially dependent order
$\alpha(\cdot)$. The system has finite propagation speed for every
\[
 c>\operatorname*{ess\,sup}_{x\in\Omega}
 \frac{1}{\sqrt{\varepsilon(x)\mu(x)}}.\qedhere
\]
\end{example}
\begin{proof}
We verify the hypotheses of \Cref{cor:bounded-memory}. For operator-valued holomorphy, we set
\[
  f_\beta(z)\coloneqq \frac{z^\beta}{1+z^\beta},
  \qquad \operatorname{Re}z>0,\quad 0<\beta<1.
\] 
As $\operatorname{Re}z^\beta>0$, one has
$|f_\beta(z)|\leq 1$, uniformly in $\beta$ and $z$.
Fix $z_0$ with $\operatorname{Re}z_0>0$ and
$0<r<\operatorname{Re}z_0$. The Taylor coefficients of
$f_\beta$ at $z_0$ satisfy, using the estimates given by Cauchy's integral formulas,
\[
  |c_n(\beta)|\leq r^{-n}
  \qquad (n\in\mathbb N_0,\ 0<\beta<1).
\]
Moreover, $\beta\mapsto c_n(\beta)$ is continuous.
For measurable $\alpha$, the coefficients $C_n\coloneqq c_n\circ\alpha$ therefore belong to
$L_\infty(\Omega)$ and satisfy
$\|C_n\|_{L_\infty}\leq r^{-n}$. Consequently,
\[
  z\longmapsto \frac{z^{\alpha(\cdot)}}
  {1+z^{\alpha(\cdot)}}
  =\sum_{n=0}^\infty C_n(z-z_0)^n
\]
defines an $L_\infty(\Omega)$-valued holomorphic function on
$\{z\in\mathbb C:\operatorname{Re}z>0\}$, since the power
series converges in $L_\infty(\Omega)$ for $|z-z_0|<r$.
It is bounded there by $1$. Thus, for
$a,\sigma\in L_\infty(\Omega)$, the corresponding Maxwell
material law $M_1$ is bounded and holomorphic with values
in $L_\infty(\Omega)^{6\times 6}$.

Moreover, since
$\|\sum_{j=1}^3n_jC_j\|=1$ for $|n|=1$,
\[
 \left\|M_0(x)^{-1/2}\left(\sum_{j=1}^3n_jA_j\right)
 M_0(x)^{-1/2}\right\|
 =\frac1{\sqrt{\varepsilon(x)\mu(x)}}.\qedhere
\]
\end{proof}

Discrete time-delay terms also fit the Maxwell material-law framework.
\begin{remark}\label{rem:maxwell-delay}
The constitutive law in \Cref{ex:fractional-maxwell} may also include
$\mathrm{e}^{-hz}b(x)1_3$ in the electric block of $M_1(z)$, where $h>0$ and
$b\in L_\infty(\Omega)$. This corresponds to the delayed term
$b(x)E(t-h,x)$, and
\[
 \sup_{\operatorname{Re}z>0}\|\mathrm{e}^{-hz}b\|_{L_\infty(\Omega)}
 \le\|b\|_{L_\infty(\Omega)}.
\]
For $b\neq0$, the inverse Laplace transform of the additional
symbol is the operator-valued measure $b(x)\,\delta_h$, supported at the
single time $h$. This is \emph{not} an element of
$L_1(\mathbb R_{\ge0};L_\infty(\Omega)^{3\times3})$:
it has a point-mass component rather than an integrable
function as its time kernel. In particular, this delay term is not
covered by the $L_1$-convolution-kernel assumption in
\cite{BlazekStolkSymes2013}. Nevertheless, $z\mapsto \mathrm{e}^{-hz}b$
is bounded and holomorphic on the right half-plane and spatially
local, so \Cref{cor:bounded-memory} applies with exactly the same propagation-speed
bound; the lower bound required for $\nu$ may change.
\end{remark}
\end{subsection}

\begin{subsection}{Elasticity with fractional stress compliance}
\label{subsec:fractional-compliance}
Fractional viscoelasticity is considered in
\cite{PicardTrostorffWaurick2015fractional}. We study here a
fractional-integral \emph{stress-compliance} relation. This is distinct
from the Kelvin--Voigt constitutive law treated in that reference.
For the here considered creep-compliance variant in linear fractional
viscoelasticity, see also \cite{MainardiSpada2011}.

Let $\Omega\subseteq\mathbb R^d$ be open, and set
$\Htopo=L_2(\Omega)^d\oplus L_{2,{\rm sym}}(\Omega)^{d\times d}$.
Let the mass density $\rho\in L_\infty(\Omega)$ be real with
$\rho\ge\rho_*>0$, and let
$C\in L_\infty(\Omega;L(\mathbb R^{d\times d}_{\rm sym}))$ be selfadjoint and uniformly
positive. We also assume
$D\in L_\infty(\Omega;L(\mathbb R^{d\times d}_{\rm sym}))$ is selfadjoint and nonnegative.
Fix $0<\alpha<1$. For displacement $u$, velocity
$v=\partial_{t,\nu}u$, and stress $T$, consider
\begin{equation}\label{eq:fractional-creep}
 \begin{cases}
 \rho\,\partial_{t,\nu}v-\operatorname{Div}T=f,\\
 \operatorname{Grad}u=C^{-1}T+D\partial_{t,\nu}^{\alpha-1}T.
 \end{cases}
\end{equation}
The operators $\operatorname{Grad}$
and $\operatorname{Div}$ are defined by
\[
\begin{array}{rcl}
\operatorname{Grad}\colon H^1(\Omega)^d \subseteq L_2(\Omega)^d
&\longrightarrow&
L_{2,\mathrm{sym}}(\Omega)^{d\times d}\\[2mm]
\phi
&\longmapsto&
\displaystyle \frac{1}{2}
\left(\partial_j\phi_k+\partial_k\phi_j\right)_{j,k\in\{1,\ldots,d\}}
\end{array}
\]
and
\[
\begin{array}{rcl}
\operatorname{Div}\colon \Htopo(\operatorname{Div},\Omega) \subseteq L_{2,\mathrm{sym}}(\Omega)^{d\times d}
&\longrightarrow&
L_{2}(\Omega)^{d}\\[2mm]
\Phi
&\longmapsto&
\displaystyle
\left(\sum_{j=1}^{d}\partial_j\Phi_{ij}\right)_{i\in\{1,\ldots,d\}},
\end{array}
\] where \[\Htopo(\operatorname{Div},\Omega)\coloneqq \left\{\Phi\in L_{2, {\rm sym}}(\Omega)^{d\times d}:  \operatorname{Div} \Phi \in L_{2}(\Omega)^{d} 
\text{ distributionally}\right\}.\]
The second equation is a fractional-integral compliance relation. 
Differentiating this relation in the framework of evolutionary equations yields
\[
 \partial_{t,\nu}\bigl(C^{-1}+D\partial_{t,\nu}^{\alpha-1}\bigr)T
 -\operatorname{Grad}v=0.
\]
With the local skew-selfadjoint realisation of
$A=\begin{psmallmatrix}0&-\operatorname{Div}\\-\operatorname{Grad}_0&0\end{psmallmatrix}$
(i.e., the standard homogeneous Dirichlet boundary conditions), this
is an evolutionary equation with
\begin{equation}\label{eq:fractional-creep-M}
 M(z)=\begin{pmatrix}
   \rho&0\\0&C^{-1}+z^{\alpha-1}D
 \end{pmatrix},\qquad
 M_0=\begin{pmatrix}\rho&0\\0&C^{-1}\end{pmatrix}.
\end{equation}

\begin{proposition}\label{prop:fractional-creep-speed}
Under the above assumptions, $(M,A)$ has finite propagation speed bounded by $c_{\rm el}$, where
\[
 c_{\rm el}\coloneqq \operatorname*{ess\,sup}_{x\in\Omega}
 \sup_{|n|=1}
 \left\|\rho(x)^{-1/2}
 \bigl(S\mapsto Sn\bigr)C(x)^{1/2}\right\|_{L(\mathcal S,\mathbb C^d)}.
\]
\end{proposition}
\begin{proof}
For each $\nu_0>0$, the function $z\mapsto z^{\alpha-1}$ is bounded and holomorphic on $\operatorname{Re}z>\nu_0$. Hence \eqref{eq:fractional-creep-M} defines a material law.
With the principal branch, $\operatorname{Re}z^\alpha\ge0$ for
$\operatorname{Re}z>0$. Hence
\[
 \operatorname{Re}(zM(z))
  = (\operatorname{Re}z)M_0
     +\begin{pmatrix}0&0\\0&(\operatorname{Re}z^\alpha)D\end{pmatrix}
  \ge (\operatorname{Re}z)M_0.
\]
The assumptions of \Cref{thm:main} hold with $c_0=0$. The principal
symbol in direction $n$ is an off-diagonal block operator, with
components $S\mapsto Sn$ and its adjoint
$v\mapsto(v\otimes n+n\otimes v)/2$.
Hence the stated value follows from \Cref{thm:main}.

It remains to check \Cref{hyp:main}(i). As $\operatorname{Grad}_0=-\operatorname{Div}^*$, the operator
$A=\begin{psmallmatrix}0&-\operatorname{Div}\\-\operatorname{Grad}_0&0\end{psmallmatrix}$ is skew-selfadjoint by \cite[Proposition~6.2.3(a)]{seifert2022evolutionary}, and
$A=\sum_{j=1}^d A_j\partial_j$ with $A_j=-\begin{psmallmatrix}0&D_j\\G_j&0\end{psmallmatrix}$, $D_jS=Se_j$ and $G_jv=\frac12(ve_j\transposed+e_jv\transposed)$. Since $D_j=G_j^*$, the matrices $A_j$ are selfadjoint. Finally, multiplication by scalar functions in $W^{1,\infty}(\Omega)$ preserves the domains of $\operatorname{Grad}_0$ and $\operatorname{Div}$, with the usual product rules, so that $A$ satisfies \eqref{eq:local-domain}.
\end{proof}

\end{subsection}

\begin{subsection}{Thermo-elasticity with dual phase-lags}  
The authors of~\cite{mukhopadhyay2011domain} investigated the domain of influence for thermo-elasticity with two-phase lags, in which thermo-elastic disturbances propagate at a finite speed. They proposed a structural threshold that serves as an upper speed bound to restrict the domain of influence, and applied an energy-type estimate to guarantee that disturbances propagate at a finite speed. In this section, we study the finite-speed propagation of the aforementioned thermo-elasticity disturbances within the evolutionary equations framework. We emphasise that the proof of finiteness of the speed of propagation does not require a sophisticated energy balance first. In fact, the evolutionary equations reformulation already contains all the necessary information. For convenience, we derive this reformulation here.

Let $\Omega\subseteq\mathbb R^d$ be a bounded open set.
Let us consider the equation of motion:
\begin{equation}
\label{eq:motion}
\mu u_{i,jj} +(\lambda+\mu) u_{j,ji} -\gamma \theta_{,i}+ F_i = \rho\ddot{u}_i
\end{equation}
and the heat conduction equation:
 \begin{equation}
 \label{eq:heat}
K\big(1+\tau_T\partial_t\big)\nabla^2\theta= \big(1+\tau_q\partial_t+\frac{1}{2}\tau_q^2\partial_t^2\big)(\rho c_v \dot{\theta}+\gamma \theta_0 \dot{u}_{i,i}-r), 
\end{equation}
on $[0,\infty)\times\Omega$.

Let us first consider the term
\begin{equation*}
\mu u_{i,jj} +(\lambda+\mu) u_{j,ji} = \mu \Delta u_i+(\lambda+\mu)\partial_i \operatorname{div} u.
\end{equation*}
Here $\lambda,\mu>0$ are the Lam\'e constants,
$\rho\in L_\infty(\Omega)$ is the mass density with
$0<\rho_*\le\rho\le\rho^*$ almost everywhere, $c_v>0$ is the
specific heat at constant strain, $K>0$ the thermal conductivity,
$\theta_0>0$ the reference temperature, $\gamma\in\mathbb R$ the
thermoelastic coupling constant, and $r$ the heat source. We take the
remaining physical parameters $c_v,K,\theta_0,\tau_T,\tau_q$ to be
positive constants, as in the homogeneous dual-phase-lag model
(see \cite[Basic equations and problem formulation]{mukhopadhyay2011domain}).
Hence, denoting the vectorial Laplacian acting on each vector component separately by $\Delta$ we get
\[\big(\mu \Delta u_i+(\lambda+\mu)\partial_i \operatorname{div} u \big)_i= \mu \Delta u+(\lambda+\mu)\operatorname{grad}\operatorname{div} u.\]
Since $\lambda,\mu>0$, we find
$0<C=C^*\in L(\mathbb{R}^{d\times d}_{\mathrm{sym}})$ strictly
positive definite (see \cite[Lemma~2.6]{Picard2015} and
\cite[Section~7.1]{seifert2022evolutionary}) such that
\[ \mu \Delta u+(\lambda+\mu)\operatorname{grad}\operatorname{div} u=\operatorname{Div} \left(C \operatorname{Grad} u\right).\]
Hence, rearranging the terms in \eqref{eq:motion}, we get
\[ \rho \partial_t^2u- \operatorname{Div} \left(C \operatorname{Grad} u\right)+\gamma  \operatorname{grad} \theta=F.\]
In this equality, we rewrite $\gamma  \operatorname{grad} \theta.$ Note that $\gamma$ is a scalar and $\theta$ is a scalar-valued function. We can also write
\[\operatorname{grad} \theta= \operatorname{Div}\operatorname{trace}^{*} \theta. \]
Therefore
\begin{equation}
\label{eq:motion-rewritten}
 \rho \partial_t^2u-\operatorname{Div}\big(C \operatorname{Grad} u-\gamma\operatorname{trace}^{*} \theta\big)=F.
\end{equation}
Next, we reformulate \eqref{eq:heat} which is
\[K \operatorname{div}\big(1+\tau_T\partial_t\big) \operatorname{grad} \theta=\big(1+\tau_q\partial_t+\frac{1}{2}\tau_q^2\partial_t^2\big)(\rho c_v \partial_t{\theta}+\gamma \theta_0 \operatorname{div}\partial_t u-r),\]
using 
\[\operatorname{div}\partial_t u=\operatorname{trace}\operatorname{grad}\partial_t u,\]
which gives
\begin{equation} \label{eq:heat-rewritten}
K \operatorname{div}\big(1+\tau_T\partial_t\big) \operatorname{grad} \theta=\big(1+\tau_q\partial_t+\frac{1}{2}\tau_q^2\partial_t^2\big)(\rho c_v \partial_t{\theta}+\gamma \theta_0 \operatorname{trace}\operatorname{grad}\partial_t u-r).\end{equation}
We divide by the thermal equation by $\theta_0$ and set
\[
 \rho_\theta\coloneqq\frac{\rho c_v}{\theta_0},\qquad
 k_\theta\coloneqq\frac{K}{\theta_0},\qquad
 r_\theta\coloneqq\frac{r}{\theta_0}.
\]
Thus $\rho_\theta$ is also bounded and uniformly positive. Dividing
by the dual-phase-lag polynomial (understood via its causal
functional calculus), we get
\[\kappa(\partial_t) \coloneqq k_\theta
(1+\tau_T\partial_t)
(1+\tau_q\partial_t+\tfrac12\tau_q^2\partial_t^2)^{-1},\] and
\begin{equation} \label{eq:te-system}
\begin{cases}
\rho \partial_t^2u-\operatorname{Div}\big(C \operatorname{Grad} u-\gamma\operatorname{trace}^{*} \theta\big)=F\\

\rho_\theta \partial_t\theta+\gamma  \operatorname{trace}\operatorname{grad}\partial_t u-\operatorname{div}\kappa(\partial_t)\operatorname{grad} \theta=r_\theta.
\end{cases}
\end{equation}
We now introduce the new unknowns 
\begin{equation}
\begin{cases}
v= \partial_t u,\\
T=C\operatorname{Grad}u-\gamma  \operatorname{trace}^{*} \theta,\\
q=\kappa(\partial_t)\operatorname{grad} \theta.
\end{cases}
\end{equation}
This yields
\begin{equation} \label{eq:te-first-order}
\begin{cases}

\rho \partial_t v-\operatorname{Div}T=F,\\

\partial_t(C^{-1}T+C^{-1}\gamma \operatorname{trace}^{*} \theta)-\operatorname{Grad}v=0,\\

\rho_\theta \partial_t\theta+\gamma  \operatorname{trace}\operatorname{grad}v-\operatorname{div}q=r_\theta,\\

\kappa(\partial_t)^{-1}q-\operatorname{grad}\theta=0.

\end{cases}
\end{equation}
Furthermore, $\operatorname{Grad}v=\partial_t(C^{-1}T+C^{-1}\gamma \operatorname{trace}^{*} \theta)$, then the third equation in \eqref{eq:te-first-order} yields
\[\rho_\theta \partial_t\theta+\partial_t(\gamma\operatorname{trace}C^{-1}T+\gamma\operatorname{trace}C^{-1}\gamma \operatorname{trace}^{*} \theta)-\operatorname{div}q=r_\theta.\]
Therefore, we obtain the system
\begin{equation} \label{eq:te-final}
\begin{cases}

\rho \partial_{t} v-\operatorname{Div}T=F,\\

\partial_{t}(C^{-1}T+C^{-1}\gamma \operatorname{trace}^{*} \theta)-\operatorname{Grad}v=0,\\

\rho_\theta \partial_{t}\theta+\partial_{t}(\gamma\operatorname{trace}C^{-1}T+\gamma\operatorname{trace}C^{-1}\gamma \operatorname{trace}^{*} \theta)-\operatorname{div}q=r_\theta,\\

\kappa(\partial_{t})^{-1}q-\operatorname{grad}\theta=0.

\end{cases}
\end{equation}
Furthermore,
\[\kappa(\partial_{t,\nu})^{-1}= \frac{1}{2k_\theta}\tau_q \partial_{t,\nu} \sigma+ \frac{1}{k_\theta}\sigma(1-\frac{1}{2}\sigma)+\frac{1}{k_\theta}\left(1-\sigma(1-\frac{1}{2}\sigma)\right)\left(1+\tau_T \partial_{t,\nu}\right)^{-1},\] $ \text{with} \ \sigma=\frac{\tau_q}{\tau_T}, \ \tau_q>0,\ \tau_T>0 \ \text{and} \ \tau_q< 2\tau_T,$ is well-defined for $\nu\ge0$. Indeed, $\operatorname{Re}\partial_{t,\nu}=\nu\ge0$ and $\tau_T>0$, therefore, $1+\tau_T \partial_{t,\nu}$ is continuously invertible.

We can then express \eqref{eq:te-final} in the following form, where we refer to the boundary conditions below
\[
\begin{aligned}
& \partial_{t,\nu} \begin{pmatrix}
\rho& 0 & 0 & 0 \\
0 & C^{-1} & C^{-1}\gamma\,\operatorname{trace}^{*} & 0 \\
0 & \operatorname{trace}\,\gamma C^{-1} &
\rho_\theta+ \gamma\operatorname{trace}\,\gamma C^{-1}\operatorname{trace}^{*} & 0 \\
0 & 0 & 0 & \partial_{t,\nu}^{-1}\kappa(\partial_{t,\nu})^{-1}
\end{pmatrix}\begin{pmatrix}
v \\
T \\
\theta \\
q
\end{pmatrix} \\
&+
\begin{pmatrix}
0 & -\operatorname{Div} & 0 & 0 \\
-\operatorname{Grad} & 0 & 0 & 0 \\
0 & 0 & 0 & -\operatorname{div} \\
0 & 0 & -\operatorname{grad} & 0
\end{pmatrix}\begin{pmatrix}
v \\
T \\
\theta \\
q
\end{pmatrix}
 =\begin{pmatrix}
F \\
0 \\   
r_\theta \\
0
\end{pmatrix}.
\end{aligned}
\]
The latter takes the form of an evolutionary equation with
\begin{equation} \label{eq:dpl-M0}
M_0 \coloneqq \begin{pmatrix}
\rho& 0 & 0 & 0 \\
0 & C^{-1} & C^{-1}\gamma\,\operatorname{trace}^{*} & 0 \\
0 & \operatorname{trace}\,\gamma C^{-1} &
\rho_\theta+ \gamma\operatorname{trace}\,\gamma C^{-1}\operatorname{trace}^{*} & 0 \\
0 & 0 & 0 & \frac{1}{2k_\theta}\tau_q\sigma
\end{pmatrix},
\end{equation}
\begin{equation}\label{eq:dpl-M1}
 M_1(z) \coloneqq \begin{pmatrix}
0 & 0 & 0 & 0 \\
0 & 0 & 0 & 0 \\
0 & 0 & 0 & 0 \\
0 & 0 & 0 & \frac{1}{k_\theta}\sigma(1-\frac{1}{2}\sigma)+\frac{1-\sigma(1-\frac{1}{2}\sigma)}{k_\theta(1+\tau_T z)}
\end{pmatrix},\qquad \operatorname{Re}z>0
\end{equation}
\ \text{and}
\begin{equation*}
A \subseteq \begin{pmatrix}
0 & -\operatorname{Div} & 0 & 0 \\
-\operatorname{Grad} & 0 & 0 & 0 \\
0 & 0 & 0 & -\operatorname{div} \\
0 & 0 & -\operatorname{grad} & 0
\end{pmatrix}.
\end{equation*}
For $A$ to become skew-selfadjoint, we need to impose some boundary conditions. For instance, one could require homogeneous Dirichlet-conditions for $v$ and $ \theta$, i.e; $v\in{\dom }(\operatorname{Grad_0})$ and $ \theta\in{\dom }(\operatorname{grad_0})$. Then, 
\begin{equation}\label{eq:te-A}
A \coloneqq \begin{pmatrix}
0 & -\operatorname{Div} & 0 & 0 \\
-\operatorname{Grad_0} & 0 & 0 & 0 \\
0 & 0 & 0 & -\operatorname{div} \\
0 & 0 & -\operatorname{grad_0} & 0
\end{pmatrix},\end{equation}
where
\[\operatorname{Grad_0}=-\operatorname{Div^*},\ \text{and}\ \operatorname{grad_0}=-\operatorname{div^*}. \]
Letting $f=\left(F,0,r_\theta,0\right)^{\transposed}$, we obtain the evolutionary equation

\begin{equation}
\label{eq:te-evolutionary}
\left(\overline{\partial_{t,\nu} M(\partial_{t,\nu})+A}\right)U=f, \quad U = \left(v, T, \theta, q\right)^{\transposed},\qquad M(z)=M_0+z^{-1}M_1(z).
\end{equation}

\begin{lemma}[{{\cite[Section~7.4]{seifert2022evolutionary}}}] \label{lem:dpl-material-law}
Let $\Htopo \coloneqq  \operatorname{L}_2(\Omega)^d\times \operatorname{L}_{2,\mathrm{sym}}(\Omega)^{d\times d}\times\operatorname{L}_2(\Omega)\times \operatorname{L}_2(\Omega)^d$ a Hilbert space and 
 \begin{align*}
M : \mathbb{C} \setminus \left\{ -\frac{1}{\tau_T},0\right\}
&\to L(\Htopo),\\
z &\mapsto\begin{pmatrix}
\rho& 0 & 0 & 0 \\
0 & C^{-1} & C^{-1}\gamma\,\operatorname{trace}^{*} & 0 \\
0 & \operatorname{trace}\,\gamma C^{-1} &
\rho_\theta+ \gamma\operatorname{trace}\,\gamma C^{-1}\operatorname{trace}^{*} & 0 \\
0 & 0 & 0 & z^{-1}\kappa(z)^{-1}
\end{pmatrix}.
\end{align*} Then $M$ defines a material law.

\end{lemma}

\begin{proof}
The only $z$-dependent entry is $z^{-1}\kappa(z)^{-1}$. Writing
$\sigma=\tau_q/\tau_T\in(0,2)$, a partial-fraction decomposition gives
\[
 z^{-1}\kappa(z)^{-1}
 =\frac{\tau_q\sigma}{2k_\theta}
 +\frac{1}{z}\left(
 \frac{\sigma(1-\sigma/2)}{k_\theta}
 +\frac{1-\sigma(1-\sigma/2)}{k_\theta(1+\tau_Tz)}
 \right).
\]
The right-hand side is holomorphic on $\mathbb C_{\operatorname{Re}>0}$.
Both coefficients in parentheses are nonnegative and sum to
$k_\theta^{-1}$. Hence, for $\operatorname{Re}z\ge\nu>0$,
\[
 \bigl|z^{-1}\kappa(z)^{-1}\bigr|
 \le\frac{\tau_q\sigma}{2k_\theta}
      +\frac{1}{k_\theta\nu}.
\]
All other blocks are bounded, proving the assertion.
\end{proof}

We state well-posedness for the thermoelastic evolutionary system.
\begin{theorem}\label{thm:te-wellposedness}
Let  $\Htopo \coloneqq  \operatorname{L}_2(\Omega)^d\times \operatorname{L}_{2,\mathrm{sym}}(\Omega)^{d\times d}\times\operatorname{L}_2(\Omega)\times \operatorname{L}_2(\Omega)^d$. Assume $\rho\in L_\infty(\Omega)$ with $\rho\ge\rho_*>0$, for some $\rho*>0$, 
$0<C=C^*\in L(\mathbb{R}^{d\times d}_{\mathrm{sym}})$, $\gamma\in\mathbb R$,
$k_\theta,\tau_q,\tau_T>0$ and $0<\sigma=\tau_q/\tau_T<2$.
With \eqref{eq:dpl-M0}, \eqref{eq:dpl-M1} and \eqref{eq:te-A}, the pair
$(M,A)$, where $M(z)=M_0+z^{-1}M_1(z)$, is evolutionary.
\end{theorem}

The proof of \Cref{thm:te-wellposedness} is based on \Cref{thm:wellposedness}. Thus, we restrict ourselves to the decisive observation in the following results. 
\begin{lemma}\label{lem:te-M0-positive}
Let $\Htopo \coloneqq  \operatorname{L}_2(\Omega)^d\times \operatorname{L}_{2,\mathrm{sym}}(\Omega)^{d\times d}\times\operatorname{L}_2(\Omega)\times \operatorname{L}_2(\Omega)^d$ be a Hilbert space, and let $M_0 \in \operatorname{L}(\Htopo)$ be given as in \eqref{eq:dpl-M0}.
Assume $\rho,\rho_\theta\in L_\infty(\Omega)$ are real-valued
and uniformly positive, $0<C=C^*\in L(\mathbb{R}^{d\times d}_{\mathrm{sym}})$,
$\gamma\in\mathbb R$, and $k_\theta,\tau_q,\sigma>0$.
Then there exists $c >0$ such that
\[
\begin{aligned}
& \operatorname{Re}\langle M_0\phi,\phi\rangle_{\Htopo }  \ge c\lVert\phi\rVert^2_{\Htopo }\quad \text{for all $\phi \in \Htopo$}.  
\end{aligned}
\]

\end{lemma}

\begin{proof}
For $U=(v,T,\theta,q)\in\Htopo$, the quadratic form of $M_0$ reads
\[
\begin{aligned}
 \langle M_0U,U\rangle_{\Htopo}
  ={}&\langle\rho v,v\rangle
     +\langle C^{-1}(T+\gamma\operatorname{trace}^*\theta),
                         T+\gamma\operatorname{trace}^*\theta\rangle\\
    &+\langle\rho_\theta\theta,\theta\rangle
     +\frac{\tau_q\sigma}{2k_\theta}\|q\|^2.
\end{aligned}
\]
All diagonal terms are uniformly positive. The triangular map
\[
 (v,T,\theta,q)\longmapsto
 (v,T+\gamma\operatorname{trace}^*\theta,\theta,q)
\]
is boundedly invertible on $\Htopo$, since $\operatorname{trace}^*$
is bounded. This yields $M_0=M_0^*\ge \eta$ for some $\eta>0$.
\end{proof}
\begin{lemma}[{{\cite[Lemma~7.4.3]{seifert2022evolutionary}}}]\label{lem:dpl-positive}
Let $0<\sigma<2$, $\tau_T,k_\theta>0$ and
\[
 m(z)=\frac{\sigma(1-\sigma/2)}{k_\theta}
 +\frac{1-\sigma(1-\sigma/2)}{k_\theta(1+\tau_Tz)}.
\]
Then, for every $\operatorname{Re}z\ge0$,
\[
 \operatorname{Re}m(z)\ge
 \frac{\sigma(1-\sigma/2)}{k_\theta}>0.
\]
\end{lemma}
\begin{proof}
For $0<\sigma<2$ both coefficients in the expression for $m(z)$
are nonnegative, and $\operatorname{Re}(1+\tau_Tz)^{-1}>0$ when
$\operatorname{Re}z\ge0$.
\end{proof}
\begin{theorem} \label{thm:te-finite-speed}
Under the hypotheses of \Cref{thm:te-wellposedness}, the evolutionary pair
$(M,A)$ from \eqref{eq:dpl-M0}, \eqref{eq:dpl-M1}, \eqref{eq:te-A} has finite speed.
\end{theorem}
\begin{remark} Applying our main theorem in the proof, we will show that \eqref{eq:cone-propagation} holds for every
\[
 c>\operatorname*{ess\,sup}_{x\in\Omega}\ \sup_{|n|=1}
 \Bigl\|M_0(x)^{-1/2}\Bigl(\sum_{j=1}^d n_jA_j\Bigr)M_0(x)^{-1/2}\Bigr\|,
\]
with $M_0$ from \eqref{eq:dpl-M0} and $A_j$ as in the proof below; see \Cref{rem:te-speed} for an explicit evaluation in a concrete example.
\end{remark}
\begin{proof}[Proof of \Cref{thm:te-finite-speed} and \Cref{thm:te-wellposedness}]
The coefficient bounds and \Cref{lem:te-M0-positive} give the uniform positivity of
$M_0$, and \Cref{lem:dpl-material-law} gives bounded holomorphy of $M$. The symbol
$M_1(z)$ in \eqref{eq:dpl-M1} is bounded and holomorphic on
$\mathbb C_{\operatorname{Re}>0}$. Moreover, \Cref{lem:dpl-positive} shows
$\operatorname{Re}M_1(z)\ge0$ on this half-plane. Consequently, for
any fixed $\nu>0$ and $\operatorname{Re}z\ge\nu$,
\[
 \operatorname{Re}\langle zM(z)U,U\rangle_{\Htopo}
 \ge\nu\langle M_0U,U\rangle_{\Htopo}
 \ge\nu c\|U\|_{\Htopo}^2.
\]
The abstract well-posedness theorem therefore applies for each
$\nu>0$. To verify the finite-propagation hypotheses, it remains
to check the spatial operator. The realisation $A$ in \eqref{eq:te-A} is skew-selfadjoint,
by $\operatorname{Grad}_0=-\operatorname{Div}^*$ and
$\operatorname{grad}_0=-\operatorname{div}^*$.
Its principal matrices are self-adjoint. Indeed,
\begin{equation*}
\begin{aligned}
& A \coloneqq \begin{pmatrix}
0 & -\operatorname{Div} & 0 & 0 \\
-\operatorname{Grad}_0 & 0 & 0 & 0 \\
0 & 0 & 0 & -\operatorname{div} \\
0 & 0 & -\operatorname{grad}_0 & 0
\end{pmatrix}\\ &\subseteq\begin{pmatrix}
0 & -\sum_{j=1}^d \partial_j D_j   & 0 & 0 \\
-\sum_{j=1}^{d}  \partial_jG_j & 0 & 0 & 0 \\
0 & 0 & 0 & -\sum_{j=1}^{d} \partial_jd_j  \\
0 & 0 & -\sum_{j=1}^d \partial_j g_j  & 0
\end{pmatrix}\\
&=\sum_{j=1}^d\partial_j A_j,
\end{aligned}
\end{equation*}
where
\begin{align*}
D_j:\R^{d\times d}_{\mathrm{sym}} &\to \R^d & G_j:\R^{d} &\to \R^{d\times d}_{\mathrm{sym}} \\
T&\mapsto Te_j & v&\mapsto \frac{1}{2}\left(ve_j^{\transposed}+e_jv^{\transposed}\right)
\end{align*} and \[\ 
A_j=-\begin{pmatrix}
0 &  D_j  & 0 & 0 \\
G_j & 0 & 0 & 0 \\
0 & 0 & 0 & d_j  \\
0 & 0 & g_j  & 0
\end{pmatrix}.
\]
Since $\langle D_jS,v\rangle=\langle S,G_jv\rangle$ for the natural
Euclidean inner products, $D_j=G_j^*$ and hence $A_j=A_j^*$. In particular, $(M,A)$ is an evolutionary pair settling the proof of \Cref{thm:te-wellposedness}.
Finally, scalar $W^{1,\infty}(\Omega)$-multiplication preserves the
domains of $\operatorname{Grad}_0$, $\operatorname{Div}$,
$\operatorname{grad}_0$ and $\operatorname{div}$ by their distributional
product rules. Hence the boundary realisation~\eqref{eq:te-A} satisfies
\eqref{eq:local-domain} as required.
Thus \Cref{cor:bounded-memory} applies and yields finite
propagation speed. 
\end{proof}

 \begin{remark}
 For illustration, we provide the matrix $A_1$ for $d=3$, where, similarly to \cite{BuchingerPorfidoWaurick2026}, we identify symmetric matrices with appropriate $6$-vectors and get
\[G_1=\begin{pmatrix}
1 & 0&0 \\
0 & \frac{1}{\sqrt{2}}&0 \\
0 & 0&\frac{1}{\sqrt{2}}\\
0&0&0\\
0&0&0\\
0&0&0
\end{pmatrix},\
G_2=\begin{pmatrix}
0 & 0&0 \\
\frac{1}{\sqrt{2}}& 0&0 \\
0 & 0&0\\
0&1&0\\
0&0&\frac{1}{\sqrt{2}}\\
0&0&0
\end{pmatrix},
\
G_3=\begin{pmatrix}
0 & 0&0 \\
0& 0&0 \\
\frac{1}{\sqrt{2}} & 0&0\\
0&0&0\\
0&\frac{1}{\sqrt{2}}&0\\
0&0&1
\end{pmatrix},
\]
\[D_1=\begin{pmatrix}
1 & 0&0&0&0&0 \\
0 & \frac{1}{\sqrt{2}}&0&0&0&0 \\
0 & 0&\frac{1}{\sqrt{2}}&0&0&0
\end{pmatrix},
\
D_2=\begin{pmatrix}
0 & \frac{1}{\sqrt{2}}&0&0&0&0 \\
0 & 0&0&1&0&0 \\
0 & 0&0&0&\frac{1}{\sqrt{2}}&0
\end{pmatrix}
\]  and \[
D_3=\begin{pmatrix}
0& 0&\frac{1}{\sqrt{2}}&0&0&0 \\
0 & 0&0&0&\frac{1}{\sqrt{2}}&0 \\
0 & 0&0&0&0&1
\end{pmatrix}.
\]
Thus,
\[
A_1=-
\begin{pmatrix}
\begin{array}{ccccccccccccc}
0&0&0&1&0&0&0&0&0&0&0&0&0\\
0&0&0&0&\frac{1}{\sqrt{2}}&0&0&0&0&0&0&0&0\\
0&0&0&0&0&\frac{1}{\sqrt{2}}&0&0&0&0&0&0&0\\
1&0&0&0&0&0&0&0&0&0&0&0&0\\
0&\frac{1}{\sqrt{2}}&0&0&0&0&0&0&0&0&0&0&0\\
0&0&\frac{1}{\sqrt{2}}&0&0&0&0&0&0&0&0&0&0\\
0&0&0&0&0&0&0&0&0&0&0&0&0\\
0&0&0&0&0&0&0&0&0&0&0&0&0\\
0&0&0&0&0&0&0&0&0&0&0&0&0\\
0&0&0&0&0&0&0&0&0&0&1&0&0\\
0&0&0&0&0&0&0&0&0&1&0&0&0\\
0&0&0&0&0&0&0&0&0&0&0&0&0\\
0&0&0&0&0&0&0&0&0&0&0&0&0
\end{array}
\end{pmatrix}
\] is self-adjoint. \qedhere
 \end{remark} 

\begin{remark}[Comparison with \cite{mukhopadhyay2011domain}: an explicit example]\label{rem:te-speed}
The admissible speed from \Cref{thm:te-finite-speed} can be strictly smaller than the one in \cite[(11)]{mukhopadhyay2011domain}. To see this, consider the homogeneous isotropic setting of \cite{mukhopadhyay2011domain}: $\rho>0$ is constant, $CS=2\mu S+\lambda\operatorname{trace}(S)1$ with $\lambda,\mu>0$, and we put
\[
 c_{\mathrm p}^2\coloneqq\frac{\lambda+2\mu}{\rho},\qquad
 \beta\coloneqq\frac{\tau_q\sigma}{2k_\theta}=\frac{\tau_q^2}{2\tau_Tk_\theta},
\]
so that $c_{\mathrm p}$ is the longitudinal wave speed ($c_1$ in \cite{mukhopadhyay2011domain}) and $\beta$ is the $q$-entry of $M_0$. We choose
\[
 \tau_q=\tau_T=\frac{2K}{\rho c_vc_{\mathrm p}^2}
 \qquad\text{and}\qquad
 \gamma^2=\frac{\rho^2c_vc_{\mathrm p}^2}{\theta_0}.
\]
Then $\sigma=1\in(0,2)$, and $\beta^{-1}=\rho_\theta c_{\mathrm p}^2$ as well as $\gamma^2=\rho\rho_\theta c_{\mathrm p}^2$. In the dimensionless variables of \cite{mukhopadhyay2011domain}, this choice reads $\tau_q=\tau_T=2$ (so their condition $\tau_q<2\tau_T$ holds) and $\varepsilon=1$. Hence \cite[(11)]{mukhopadhyay2011domain} requires $v\ge\max\{2,1+\varepsilon,2\tau_T/\tau_q^2\}=\max\{2,2,1\}=2$, that is, a propagation speed of $2c_{\mathrm p}$ in physical units.

We now compute the bound of \Cref{thm:te-finite-speed} for these parameters. Fix $n\in\R^d$ with $|n|=1$ and put $A_n\coloneqq\sum_{j=1}^dn_jA_j$. As $M_0^{-1/2}A_nM_0^{-1/2}$ is self-adjoint, its norm is the largest $|s|$ for which $A_nU=sM_0U$ has a nontrivial solution $U=(v,T,\theta,q)$. For $s\neq0$ and $S\coloneqq C^{-1}(T+\gamma\operatorname{trace}^*\theta)$, this eigenvalue problem reads
\[
 -Tn=s\rho v,\quad
 -\tfrac12\bigl(vn\transposed+nv\transposed\bigr)=sS,\quad
 -n\cdot q=s\bigl(\rho_\theta\theta+\gamma\operatorname{trace}S\bigr),\quad
 -\theta n=s\beta q.
\]
For $v\perp n$, one finds $\operatorname{trace}S=0$, $\theta=0$, $q=0$ and $s^2=\mu/\rho<c_{\mathrm p}^2$. Indeed, the second
of the four eigenvalue equations provides
\[
 S=-\frac{1}{2s}(vn\transposed+nv\transposed),\qquad
 \operatorname{trace}S=-\frac{v\cdot n}{s}=0,\qquad
 Sn=-\frac{1}{2s}v.
\]
By isotropy and the definition of $S$,
\[
 T=CS-\gamma\theta 1
  =2\mu S+\lambda\operatorname{trace}(S)1-\gamma\theta1
  =2\mu S-\gamma\theta1.
\]
Inserting this into the first eigenvalue equation yields
\[
 \frac{\mu}{s}v+\gamma\theta n=s\rho v.
\]
Taking the scalar product with $n$ and using $v\perp n$ implies
$\gamma\theta=0$. Since $\gamma\neq0$ for the present parameter
choice, $\theta=0$; the fourth equation then implies $q=0$.
The remaining identity $\mu v/s=s\rho v$, with $v\neq0$, reads
\[
 s^2=\frac{\mu}{\rho}<c_{\mathrm p}^2.
\]
 For $v=wn$, the second equation gives $S=-\frac ws nn\transposed$, so that the first one yields $\gamma\theta=\frac ws\bigl(\rho s^2-(\lambda+2\mu)\bigr)$, and the last two give $\theta\bigl(\beta^{-1}-\rho_\theta s^2\bigr)=-s\gamma w$. A nontrivial solution therefore exists if and only if
\[
 \bigl(\rho s^2-(\lambda+2\mu)\bigr)\bigl(\rho_\theta s^2-\beta^{-1}\bigr)=\gamma^2s^2.
\]
Dividing by $\rho\rho_\theta c_{\mathrm p}^4$ and inserting the above choice of parameters, we obtain for $\xi\coloneqq s^2/c_{\mathrm p}^2$
\[
 (\xi-1)^2=\xi,\qquad\text{i.e.}\qquad \xi^2-3\xi+1=0,
\]
whose larger root is $\xi=\frac{3+\sqrt5}{2}=\phi^2$, where $\phi=\frac{1+\sqrt5}{2}$. Hence, for every $n$,
\[
 \bigl\|M_0^{-1/2}A_nM_0^{-1/2}\bigr\|=\phi\,c_{\mathrm p},
\]
and \Cref{thm:te-finite-speed} yields finite propagation speed for every $c>\phi\,c_{\mathrm p}\approx1.618\,c_{\mathrm p}$, whereas \cite[(11)]{mukhopadhyay2011domain} gives $2c_{\mathrm p}$. Moreover, the analysis in \cite{mukhopadhyay2011domain} assumes that the displacement is the gradient of a scalar potential (see~\cite[(4)--(7)]{mukhopadhyay2011domain}), whereas our result allows general displacement fields and bounded measurable, uniformly positive densities, and is formulated for weighted $L_{2,\nu}$-solutions.
\end{remark}
\end{subsection} 

\end{section}

\section*{Declaration for the use of AI tools}
During the final phase of preparing this manuscript, ChatGPT (including
GPT-5.6 and GPT-6) and Claude Opus 5.5 were used to assist with
proofreading, improving the presentation, checking consistency, and
exploring alternative proof strategies. The authors independently
verified the mathematical arguments, all claims and results.


\begin{thebibliography}{99}

\bibitem{CombesThomas1973}
J.~M. Combes and L.~Thomas.
\newblock \emph{Asymptotic behaviour of eigenfunctions for multiparticle Schr\"odinger operators}.
\newblock Communications in Mathematical Physics, 34:251--270, 1973.

\bibitem{CoulhonSikora2008}
T.~Coulhon and A.~Sikora.
\newblock \emph{Gaussian heat kernel upper bounds via the Phragm\'{e}n--Lindel\"of theorem}.
\newblock Proceedings of the London Mathematical Society, 96(2):507--544, 2008.



\bibitem{antontsev2016finite}
S.~N. Antontsev and J.~I. D{\'\i}az.
\newblock {\em Finite speed of propagation and waiting time for local solutions of degenerate equations in viscoelastic media or heat flows with memory}.
\newblock {\em J. Elliptic Parabolic Equ.}, 2(1):207--216, 2016.

\bibitem{BuchingerSkrepekWaurick2024}
A.~Buchinger, N.~Skrepek, and M.~Waurick.
\newblock {\em Weak Operator Continuity for Evolutionary Equations}.
\newblock Pure and Applied Functional Analysis, 9(4):963--990, 2024.

\bibitem{BuchingerPorfidoWaurick2026}
A. Buchinger, M. Porfido, and M. Waurick,
\emph{The Free Lunch Theorem of Homogenisation},
arXiv:2606.04498 (2026).
\url{https://arxiv.org/abs/2606.04498}.

\bibitem{cattaneo1958forme}
C.~Cattaneo.
\newblock {\em Sur une forme de l'équation de la chaleur éliminant le paradoxe d'une propagation instantanée}.
\newblock Compt. Rendu, 247:431--433, 1958.

\bibitem{cheng2023dynamics}
C.~Cheng.
\newblock {\em Dynamics and spreading speeds of a nonlocal diffusion model with advection and free boundaries}.
\newblock arXiv preprint arXiv:2309.06185, 2023.

\bibitem{Christov2009}
C.~I. Christov.
\newblock {\em On frame indifferent formulation of the Maxwell--Cattaneo model of finite-speed heat conduction}.
\newblock Mechanics Research Communications, 36(4):481--486, 2009.

\bibitem{FinnWheeler1972}
J.~M. Finn and L.~T. Wheeler.
\newblock {\em Wave propagation aspects of the generalized theory of heat conduction}.
\newblock Zeitschrift für Angewandte Mathematik und Physik, 23:927--940, 1972.

\bibitem{Friedrichs1954}
K.~O. Friedrichs.
\newblock {\em Symmetric hyperbolic linear differential equations}.
\newblock Communications on Pure and Applied Mathematics, 7(2):345--392, 1954.

\bibitem{Friedrichs1958}
K.~O. Friedrichs.
\newblock {\em Symmetric positive linear differential equations}.
\newblock Communications on Pure and Applied Mathematics, 11(3):333--418, 1958.

\bibitem{Goldstein1951}
S.~Goldstein.
\newblock {\em On diffusion by discontinuous movements, and on the telegraph equation}.
\newblock The Quarterly Journal of Mechanics and Applied Mathematics, 4(2):129--156, 1951.

\bibitem{grimmer1984energy}
R.~Grimmer and M.~Zeman.
\newblock {\em Energy estimates for symmetric hyperbolic integro-differential equations}.
\newblock North-Holland Mathematics Studies, 92:241--248, 1984.

\bibitem{GurtinPipkin1968}
M.~E. Gurtin and A.~C. Pipkin.
\newblock {\em A general theory of heat conduction with finite wave speeds}.
\newblock Archive for Rational Mechanics and Analysis, 31(2):113--126, 1968.

\bibitem{ignaczak2010thermoelasticity}
J.~Ignaczak and M.~Ostoja-Starzewski.
\newblock {\em Thermoelasticity with Finite Wave Speeds}.
\newblock Oxford University Press, Oxford, 2010.

\bibitem{Keller2004}
J.~B. Keller.
\newblock {\em Diffusion at finite speed and random walks}.
\newblock Proceedings of the National Academy of Sciences, 101(5):1120--1122, 2004.

\bibitem{BlazekStolkSymes2008}
K.~D. Blazek, C.~Stolk, and W.~W. Symes.
\newblock \emph{Well-posedness of initial/boundary value problems for
hyperbolic integro-differential systems with nonsmooth coefficients}.
\newblock Preprint, 2008.
\newblock \url{https://repository.rice.edu/server/api/core/bitstreams/7db2cca2-46ed-41f4-a8db-7306dd56be52/content}.

\bibitem{BlazekStolkSymes2013}
K.~D. Blazek, C.~Stolk, and W.~W. Symes.
\newblock \emph{A mathematical framework for inverse wave problems in
heterogeneous media}.
\newblock Inverse Problems, 29:065001, 2013.
\newblock arXiv:1211.4656. \url{https://arxiv.org/abs/1211.4656}.

\bibitem{FinsterMurroSchmid2025}
F.~Finster, S.~Murro, and G.~Schmid.
\newblock \emph{The Cauchy problem for symmetric hyperbolic systems with
nonlocal potentials}.
\newblock arXiv:2507.05004, 2025. \url{https://arxiv.org/abs/2507.05004}.

\bibitem{MainardiSpada2011}
F.~Mainardi and G.~Spada.
\newblock \emph{Creep, relaxation and viscosity properties for basic fractional models in rheology}.
\newblock The European Physical Journal Special Topics, 193:133--160, 2011.
\newblock \url{https://arxiv.org/abs/1110.3400}.

\bibitem{McIntoshMorris2013}
A.~McIntosh and A.~J. Morris.
\newblock {\em Finite propagation speed for first order systems and Huygens' principle for hyperbolic equations}.
\newblock Proceedings of the American Mathematical Society, 141(10):3515--3527, 2013.

\bibitem{McLaughlinYoon2017}
J.~R. McLaughlin and J.-R. Yoon.
\newblock {\em Finite propagation speed of waves in anisotropic viscoelastic media}.
\newblock SIAM Journal on Applied Mathematics, 77(6):1921--1936, 2017.

\bibitem{moussaoui2020speed}
A.~Moussaoui and V.~Volpert.
\newblock {\em Speed of wave propagation for a nonlocal reaction--diffusion equation}.
\newblock {\em Appl. Anal.}, 99(13):2307--2321, 2020.

\bibitem{mukhopadhyay2011domain}
S.~Mukhopadhyay, S.~Kothari, and R.~Kumar.
\newblock {\em A domain of influence theorem for thermoelasticity with dual phase-lags}.
\newblock Journal of Thermal Stresses, 34(9):923--933, 2011.

\bibitem{PicardTrostorffWaurick2015fractional}
R.~Picard, S.~Trostorff, and M.~Waurick.
\newblock {\em On evolutionary equations with material laws containing fractional integrals}.
\newblock Mathematical Methods in the Applied Sciences, 38(15):3141--3154, 2015.
\newblock doi:10.1002/mma.3286.

\bibitem{Picard2015}
R.~Picard, S.~Seidler, S.~Trostorff, and M.~Waurick.
\newblock {\em On abstract $\mathrm{grad}$--$\mathrm{div}$ systems}.
\newblock J. Differential Equations, 260(6):4888--4917, 2016.
\newblock doi:10.1016/j.jde.2015.11.033.

\bibitem{picard2009structural}
R.~Picard.
\newblock {\em A structural observation for linear material laws in classical mathematical physics}.
\newblock Mathematical Methods in the Applied Sciences, 32(14):1768--1803, 2009.
\newblock Wiley.

\bibitem{PicardMcGhee2011}
R.~Picard and D.~McGhee.
\newblock {\em Partial Differential Equations: A Unified Hilbert Space Approach}.
\newblock De Gruyter Expositions in Mathematics, volume 55.
\newblock De Gruyter, Berlin, 2011.

\bibitem{picard2017maximal}
R.~Picard, S.~Trostorff, and M.~Waurick.
\newblock {\em On maximal regularity for a class of evolutionary equations}.
\newblock Journal of Mathematical Analysis and Applications, 449(2):1368--1381, 2017.

\bibitem{seifert2022evolutionary}
C.~Seifert, S.~Trostorff, and M.~Waurick.
\newblock {\em Evolutionary Equations: Picard's Theorem for Partial Differential Equations, and Applications}.
\newblock Operator Theory: Advances and Applications, volume 287.
\newblock Birkh{\"a}user, Cham, 2022.

\bibitem{segatti2026alpha}
A.~Segatti and R.~Taranets.
\newblock {\em On $\alpha$-entropy solutions of a nonlocal thin film equation: existence and finite speed of propagation}.
\newblock arXiv preprint arXiv:2601.19808, 2026.

\bibitem{stana2016finite}
D.~Stan, F.~del~Teso, and J.~L.~V{\'a}zquez.
\newblock {\em Finite and infinite speed of propagation for porous medium equations with nonlocal pressure}.
\newblock {\em J. Differ. Equations}, 260(2):1154--1199, 2016.

\bibitem{Trostorff2015}
S.~Trostorff.
\newblock {\em Exponential Stability for Second Order Evolutionary Problems}.
\newblock Journal of Mathematical Analysis and Applications, 429(2):1007--1032, 2015.

\bibitem{TrostorffWaurick2021}
S.~Trostorff and M.~Waurick.
\newblock {\em Maximal Regularity for Non-Autonomous Evolutionary Equations}.
\newblock Integral Equations and Operator Theory, 93 (2021), no.~2, Paper No.~30, 37 pp.
\newblock arXiv:2006.16696.

\end{thebibliography}
\end{document}